\documentclass[11pt,4paper]{amsart}

\usepackage{amssymb}
\usepackage{bbm}
\usepackage{amsmath}
\usepackage{amsthm}
\usepackage{enumitem}
\usepackage{amsfonts}
\usepackage{hyperref}
\usepackage{cleveref}
\usepackage{comment}
\usepackage{xcolor}
\usepackage{cancel} 

\usepackage{graphicx}

\usepackage{layout}

\usepackage{bbm}
\usepackage{amssymb}
\usepackage{mathtools}
\usepackage{cleveref}

\makeatletter
\@namedef{subjclassname@2020}{%
  \textup{2020} Mathematics Subject Classification}
\makeatother

\usepackage[T1]{fontenc}

\usepackage{geometry}
\usepackage{setspace}
\newtheorem{thm}{Theorem}[section]
\newtheorem{crl}[thm]{Corollary}
\newtheorem{lem}[thm]{Lemma}
\newtheorem{prop}[thm]{Proposition}

\newtheorem{example}[thm]{Example}

\theoremstyle{definition}
\newtheorem{dfn}[thm]{Definition}
\newtheorem{rem}[thm]{Remark}
\newtheorem{notation}[thm]{Notation}
\newtheorem{fact}[thm]{Fact}

\newtheorem*{convention}{Convention}
\newtheorem*{assumption}{Assumption}

\newtheorem*{qst*}{Question}

\numberwithin{equation}{section}

\DeclareMathOperator{\opaxsem}{SEM}
\DeclareMathOperator{\opaxisem}{ISEM}

\DeclareMathOperator{\cl}{cl}
\DeclareMathOperator{\nat}{nat}

\DeclareMathOperator{\iMe}{i}
\DeclareMathOperator{\wMe}{w}
\DeclareMathOperator{\diag}{Diag}
\DeclareMathOperator{\identity}{id}
\DeclareMathOperator{\acl}{acl}

\DeclareMathOperator{\Mod}{Mod}

\DeclareMathOperator{\nsop}{NSOP}
\DeclareMathOperator{\sop}{SOP}
\DeclareMathOperator{\trp}{TP}
\DeclareMathOperator{\ntp}{NTP}

\DeclareMathOperator{\tr}{Tr}
\DeclareMathOperator{\fan}{Fan}
\DeclareMathOperator{\spath}{SPath}

\newcommand{\kpluszero}{\mathcal{K}_{\scalebox{.5}{0}}} 
\newcommand{\kpluszerobar}{\overline{\mathcal{K}}_{\scalebox{.5}{0}}}

\newcommand{\kmuinfinity}{\mathcal{K}_{\rule[.8ex]{0ex}{0ex}\scalebox{.6}{$ \mu $}}^{\rule[-.4ex]{0ex}{0ex}\scalebox{.6}{$ \infty $}}}

\newcommand{\tmuinfinity}{T_{\rule[.8ex]{0ex}{0ex}\scalebox{.6}{$ \mu $}}^{\rule[-.4ex]{0ex}{0ex}\scalebox{.6}{$ \infty $}}}

\newcommand{\tplus}{T_{\rule[.8ex]{0ex}{0ex}\scalebox{.6}{$ \nat $}}^{\rule[-.4ex]{0ex}{0ex}\scalebox{.6}{$ \infty $}}}

\newcommand{\tmu}{\mathbb{T}_{\rule[.8ex]{0ex}{0ex}\scalebox{.6}{$ \mu $}}}

\newcommand{\uccext}[1][M]{\langle #1\rangle_{\rule[.8ex]{0ex}{0ex}\scalebox{.6}{$ \mu $}}^{\rule[-.4ex]{0ex}{0ex}\scalebox{.6}{$ \infty $}}}

\newcommand*{\intrext}[2][A]{#1<_{\iMe}#2}
\newcommand*{\wintrext}[2][A]{#1\leq_{\iMe}#2}
\newcommand{\axsem}{$ \opaxsem $\hspace*{3pt}}
\newcommand{\axisem}{$ \opaxisem $\hspace*{3pt}}
\newcommand*{\card}[1][M]{\vert #1\vert}
\newcommand{\subsetfinite}{\subseteq_{\omega}}
\renewcommand{\restriction}{\mathord{\upharpoonright}}
\newcommand{\wcl}{\cl^{\wMe}}
\newcommand*{\freejoin}[3]{#1\otimes_{\hspace*{-3pt}\vspace*{8pt}#2}#3}

\begin{document}


\title[A Tame Generic Structure with Non-Algebraic Geometric Closure]{A Tame Generic Structure with Non-Algebraic Geometric Closure}

\author[S. Jalili]{Somaye Jalili}
\address{School of Mathematics\\ Institute for Research in Fundamental Sciences (IPM)\\Tehran \\
	19395-5746, Iran}
\email{s.jalili@ipm.ir}

\author[M. Pourmahdian]{Massoud Pourmahdian}
\address{School of Mathematics\\ Institute for Research in Fundamental Sciences (IPM)\\Tehran \\
	19395-5746, Iran }
\email{pourmahd@ipm.ir}

\author[A. N. Valizadeh]{Ali N. Valizadeh}
\address{Department of Mathematics and Statistics\\ University of Isfahan\\Isfahan\\
	81746-73441, Iran}
\email{a.valizadeh@mcs.ui.ac.ir}

\date{\today}

\begin{abstract}
		By providing a procedure to apply Hrushovski's amalgamation method to the setting of classes of infinite structures, we introduce the notion of \textit{paracollapsed} structures. We show that this approach provides existentially closed generic structures in which the geometric closure is not included in the algebraic closure while the resulting theory is decidable. We show that paracollapsed structures have the strict order property and $\text{TP}_2$.
	    
%
%
%
\end{abstract}

\subjclass[2020]{03C45}

\keywords{Unstable theories, Strong order property, strict order property (SOP), generic structures, Hrushovski constructions, Fra\"iss\'e-Hrushovski limits, amalgamation classes.}

\maketitle

\section{Introduction}

Generic structures obtained from Hrushovski's amalgamation method can be divided into two different classes according to the way that the \textit{geometric closure} interacts with the algebraic closure. As will further be elaborated below, generics in which geometric closure in included in the algebraic closure possess tamer theories comparing to the opposite situation--where building blocks of geometric closure are not algebraic--that leads to theories which are too wild by being undecidable and having the strict order property. By resolving the undecidability issue, this work forms the first step towards finding Hrushovski constructions in which geometric closure is neither trivial nor included in the algebraic closure, but the resulting theory is computationally tame, i.e. decidable.

Geometric closure in these structures is derived from a well-behaved predimension function that assigns to each finite structure, $ A $, a real number that indicates a certain measure of sparsity for the structure. All finite structures $ B\supseteq A $ with non-positive relative predimension over $ A $ are contained in the geometric closure of $ A $. In other words, the geometric closure of $ A $ accumulates all finite superstructures $ B $ which increase the overall density of the accumulated structure (see \Cref{secPrelim}).

Although other factors are involved, a wide range of model-theoretic properties of generic structures obtained by this method is influenced by the interaction of geometric closure with the algebraic closure.‌ As said earlier, tamer generic structures are obtained in situations where geometric closure is small in the sense that it is contained in the algebraic closure: Hrushovski's new strongly minimal set (\cite{Hrushovski-NewStrongly}), stable and $\omega$-stable generics studied in \cite{BaldwinShi-StableGen} and \cite{ValiPourmahd-SSS}, simple structures obtained using a collapsed function (\cite{Pourmahd-SimpleGen}), and $ \nsop_4 $ generics obtained from a control function (\cite{EvansWong-SomeRemarksonGen}). From a different perspective, it is shown in \cite{ValiPourmahd-PSH} that a transition from binary to ternary relations can dramatically change the pseudofiniteness of stable generic structures.

The flexible nature of this construction method, first introduced in \cite{Hrushovski-NewStrongly}, gives it the potentiality of  being deployed into diverse model-theoretic settings. Beside the theories mentioned above, worthy of note are the recent results appeared in \cite{Barbina&Casanovas-STS} where it is proved that Fra\"iss\'e limit of Steiner triple systems is $\nsop_1$ and has $\trp_2$. Also, it is shown in \cite{Baldwin&Paolini-Strongly_minial_Steiner_systemsI_Existence} and \cite{Baldwin-STS} that applying Hrushovski's method to Steiner triple systems leads to finding a range of new kind of strongly minimal counterexamples for Zilber's trichotomy conjecture. 

%

As mentioned earlier, in situations where the geometric closure exceeds the bounds of algebraic closure, this method yields extremely wild theories which have the strict order property (\cite{EvansWong-SomeRemarksonGen}), and are essentially undecidable in the sense that even an incomplete fragment of them interprets Robinson arithmetic (\cite{BrodyLaskowski-OnRationalLimits}). 

\subsection{Two ways of defining geometric closure.} As evident, the notion of geometric closure is primarily based on the notion of a substructure \textit{being closed} in a larger one. And, as will be elaborated in \Cref{subsecHru}, there are two general ways of defining this notion in the context of Hrushovski's amalgamation methods derived from a predimension; we use the symbols $ A\leq B $ and $ A<B $ to denote them, which are here referred to as ``A is weakly closed in B'' and ``A is closed in B'', respectively. For a given finite substructure $A$, all finite structures $B$ containing $A$, in which $A$ is not (weakly) closed constitute the main building blocks of the (weak) geometric closure of $A$.  In cases where $A$ is not weakly closed in $B$, that is $A\not\leq B$, the structure $B$ contains a substructure whose relative predimension over $A$ is negative. This fact ensures that the weak geometric closure is contained in the algebraic closure (see \Cref{factUpperBound}).

However, the geometric closure may involve other kinds of building blocks to which we refer as \textit{zero minimal extensions} (see \Cref{dfnSS}). For a finite structure $A$, these minimal extensions concern finite structures $B\supseteq A$ in which $A$ is weakly closed but not closed, that is $A\leq B$ and $A\not<B$; such extensions will have a zero relative predimension over $A$. This fact, generally, allows zero minimal extensions to find infinitely many realizations over $A$ in the generic structure, except if they are controlled using a criteria forced on the class of finite structures with which we start the machinery. From the perspective that concerns us here, what distinguishes different generic structures, first of all is whether or not the zero minimal extensions are allowed to participate in geometric closure. Secondly, and in cases where zero minimal extensions constitute part of the geometric closure, the question would be if they are controlled in a way or else they left to the excess.

Generics which excludes zero minimal extensions in their geometric closure all ly in the realm of stable structures mainly due to the fact that the number of types is controlled by the upper bound intrinsic to negative minimal extensions (see Theorem 3.34 in \cite{BaldwinShi-StableGen}, and Lemma 1.30 of \cite{BaldwinShelah-Randomness})). 


 In the case that zero minimal extensions are taking part in the geometric closure, we are faced with an interesting dichotomy: EITHER the number of realizations of zero minimal extensions is forced to be finite, and hence, the geometric closure is contained in the algebraic closure, which leads to having generic structures which are decidable, and either simple or $ \nsop_4 $ (\cite{Pourmahd-SimpleGen}, \cite{EvansWong-SomeRemarksonGen}). OR, the construction allows non-algebraic zero minimal extensions to exist, which yields generic structures which are essentially undecidable and have the strict order property (\cite{EvansWong-SomeRemarksonGen}, \cite{BrodyLaskowski-OnRationalLimits}). In a generic structure, the geometric closure of a set determines its type, and as long as the construction allows the geometric closure to grow unconditionally, the types over sets become more and more diverse. This primarily accounts for why the latter mentioned structures and theories are so wild.



The notion of a \textit{paracollapsed} structure is an attempt towards merging the available approaches by a method that can roughly be seen as taking an average. However, the combinatorial configurations stemming from non-algebraic zero minimal extensions turns out to weigh much more than what might initially be expected. We prove that the limit paracollapsed structures, though decidable and highly weakened through the process, still produce configurations which are able to witness the strict order property or $\trp_2$. Towards introducing the limit models, we axiomatize the class of existentially closed models of a suitable  theory rather than following exactly the amalgamations routine of Hrushovski constructions; further details will follow, see also \Cref{rmk-Hrushovski-amalgamation}.

The next step in progress for this study is to find further adjustments which can provide non-simple theories that do not witness strong order properties of certain kinds. To better grasp the consequences of the presence of non-algebraic minimal extensions in paracollapsed generics, we also include a subsection that provides examples of formulas witnessing different model-theoretic properties like $\trp_2 $ and $ \sop_{2}$. Also, in \Cref{ex-sopn}, we include an example that introduces a formula witnessing $\sop_{n}$ for all $n\in\omega$. Below, we outline the procedure leading to defining the theory $\tmu$ that axiomatizes the class of existentially closed paracollapsed structures:

Starting by a certain natural-valued function $ \mu $, we introduce a class of infinite structures that are called $ \mu $-\textit{paracollapsed}, or paracollapsed in short. This class is axiomatized by the theory $ \tmuinfinity $. Based on this class of structures, we define a theory, $ \tmu $, whose axioms are obtained from a suitable adaptation of what is called \textit{semigenericity}, which will come along with a new axiom scheme--called \textit{intrinsic semigenericity}--that extends the scope of the axiom of semigenericty to configurations that constitute the non-algebraic part of pregeometric closures.  We will show that $ \tmu $ is the model completion of an inductive theory which is denoted by $ \tplus $. The latter intermediary theory is an expansion by definition of $ \tmuinfinity $.

\subsection{Organization of the Paper.}
 In Section \ref{secPrelim}, we first provide a concise account on Hrushovski constructions that are built using \textit{predimension} functions. This section also contains improvements of some of the well-known notions already available in the literature. We proceed in subsection \ref{subsecNat} by recalling an expansion by definition of the language that will be denoted by $\mathcal{L}^{+}$ which will be associated with a theory called $ T_{\nat} $. The latter theory is intended to guarantee a naturally expected interpretation for the new augmented predicates of $ \mathcal{L}^{+} $.

In Section \ref{secTmu}, we define $ \tmu $ that is built upon the notion of a $ \mu $-paracollapsed structure. One of the main goals of this section is to prove the consistency of $ \tmu $, addressed in subsection \ref{subsecCons}. In fact, via two main steps in Lemmas \ref{lmaUCCExtension} and \ref{lmaExCloedIsUCCGeneric}, we show that any existentially closed model of the inductive $ \mathcal{L}^{+} $-theory $ T_{\nat}\cup\tmuinfinity $ is a model of $ \tmu $. We continue in subsection \ref{subsecQE} by proving that $ \tmu $ admits a quantifier elimination down to the predicates available in $ \mathcal{L}^{+} $ (Theorem \ref{thmQE}). This, using the results of subsection \ref{subsecCons}, ultimately implies that $ \tmu $ fully axiomatizes the class of existentially closed models of $ T_{\nat}\cup\tmuinfinity $. 

The last section contains our theorems showing that $\tmu$ has the strict order property and $\trp_2$. This section contains several examples of formulas witnessing IP, $\sop_2$, and a formula that has $\sop_{n}$ for each $n\geq3$. Also, Lemma \ref{lmaIndisc} in this section distinguishes a key property of indiscernible sequences in the context of Hrushovski constructions built upon a predimension function.

\vspace*{.3cm}
\section{Preliminaries}\label{secPrelim}
The necessary preliminaries on Hrushovski constructions will be provided in subsection \ref{subsecHru}. In subsection \ref{subsecNat}, we will review a \textit{natural} expansion by definition of the language that forms the initial step towards eliminating the quantifiers in $ \tmu $.

Throughout, we will be working in a language $ \mathcal{L} $ that consists of only a ternary relation $ R $ whose interpretation $ R^{A} $ is always irreflexive and symmetric; that is, $ R^{A} $ holds only for tuples of distinct elements and whenever it holds for an ordered tuple $ (x,y,z) $, then any other permutation of $ x,y $ and $ z $ will also satisfy $ R^{A}. $ Hence, by abuse of notation, we denote by $ \card[R^{A}] $ the number of $ 3 $-element sets rather than the ordered triples satisfying $ R $ in $ A $.

\begin{notation}
	By $ A, B, C, \ldots $ we denote finite $ \mathcal{L} $-structures while $ M,N, \ldots $ represent arbitrary structures that might be infinite. By $ A\subsetfinite M, $ we mean that $ A $ is a finite substructure of $ M $.
\end{notation}
\begin{dfn}
	For structures $ M_0,M_1 $ and $ M_2 $ with $ M_1\cap M_2=M_0 $, the \textit{free join} or the \textit{free amalgamation} of $ M_1 $ with $ M_2 $ over $ M_0 $ is denoted by $ \freejoin{M_1}{M_0}{M_2} $ and is defined as the $ \mathcal{L} $-structure with universe $ M_1\cup M_2 $ whose set of relations consists only of $ R^{M_1}\cup R^{M_2} $. In other words, no relation having components both in $ M_1\backslash M_0 $ and $ M_2\backslash M_0 $ does exist in the free join.
\end{dfn}

\vspace*{.3cm}
\subsection{Hrushovski Constructions Using a Predimension}\label{subsecHru}\hfil

An integer-valued function $ \delta, $ called  the \textit{predimension} function, is defined over all finite structures. This function measures in a way the sparsity of a finite $ \mathcal{L} $-structure; having more relations in a structure leads to smaller assigned values by $ \delta $. Using this predimension function, the class of all finite $ \mathcal{L} $-structures can be equipped by two binary relations $ \leq $ and $ < $ each strengthening the ordinary notion of a substructure.

For any finite $ \mathcal{L} $-structure $ A $, let
\[ \delta(A):=\card[A]-\card[R^{A}]. \]
For $A$ and $B$ substructures of an $ \mathcal{L} $-structure $ M $, the \textit{relative predimension} of $ B $ over $ A $ is defined as
\[ \delta(B/A):=\delta(AB)-\delta(A), \]
where $ AB $ denotes the structure that is induced by $ M $ on $ A\cup B. $ Note that the relative predimension is defined regardless of whether $A$ and $B$ intersect one another or not.

\begin{dfn}\label{dfnClosed}
	Let $ A\subseteq B $ be finite.
	
	\begin{itemize}
		\item[(i)] $ A $ is said to be \textit{weakly  closed}, or \textit{w-closed}, in $ B $ if for every structure $ C $ with $ A\subseteq C \subseteq B $ we have that $ \delta(C/A)\geq0. $ This is denoted by $ A\leq B. $ 
		
		\item[(ii)] $ A $ is called \textit{strictly closed}, or \textit{closed}, in $ B $ if for every structure $ C $ with $ A\subsetneq C \subseteq B $ we have that $ \delta(C/A)>0. $ In notations we write $ A< B. $ 
		
		\item[(iii)] We say that $ A $ is $ n $-\textit{weakly closed} in $ B, $ denoted by $ A\leq_{n} B, $ if for every structure $ C $ with $ A\subseteq C \subseteq B $ and $ \card[C\backslash A]\leq n $ we have that $ A\leq C. $
		
		\item[(iv)] For $ A\subsetfinite M, $ the substructure $ A $ is said to be closed/w-closed in $ M $ if $ A $ is closed/w-closed in any finite substructure $ B\subseteq M $ containing $ A. $ These notions are  respectively denoted by $ A<M $ and $ A\leq M. $ Similarly, we can define $ A $ being $ n $-weakly closed in $ M $ which is written as $ A\leq_{n}M $.
		
	\end{itemize}
\end{dfn}
Given $ A\subsetfinite M, $ all intermediate finite structures $ A\subseteq C\subsetfinite M $ that prevent $ A $ from being closed/w-closed in $ M $ are accumulated in a set that is called the \textit{closure/w-closure} of $ A $ in $ M. $ The following definitions and remarks provide a uniform way for describing the building blocks of closure/w-closure of a set in the present context.
\begin{dfn}\label{dfnSS}
	Let $ A\subseteq B $.
	\begin{itemize}
		\item[(i)] 
		$ B $ is called a \textit{zero minimal extension} of $ A $, 
		if $ A\not<B $ but $ A\leq B $ and for every $ C $ with $ A\subseteq C\subsetneq B $ we have that $ A<C. $
		
		
		\item[(ii)]
		We say that $ B $ is a \textit{negative minimal extension} of $ A $, 
		if $ A\not\leq B $ and for every $ C $ with $ A\subseteq C\subsetneq B $ we have that $ A<C. $ 
		
	\end{itemize}
	
\end{dfn}

\begin{rem}
	The relative predimension of a zero minimal extension over the base set is zero while for $ B $ a negative minimal extension of $ A $ we have that $ \delta(B/A)<0. $
\end{rem}

\begin{dfn}\label{dfnIntrExt}
	Let $ A\subseteq B $ be finite. $ B $ is said to be an \textit{intrinsic extension} of $ A, $ denoted by $ \intrext{B}, $ if it is the union of a chain of structures like
	\[ A=B_{0}\subseteq B_{1}\subseteq\cdots\subseteq B_{n}=B, \]
	where for each $ i\geq 1, $ the structure $ B_i $ is a minimal extension of $ B_{i-1}. $ An intrinsic extension is called a \textit{weak intrinsic extension}, in notations $ \wintrext{B} $, if the above chain ends in a negative minimal extension.
	
\end{dfn}
The facts below are easily followed from the definitions.
\begin{fact}\label{factIntrExt}
	Let $ A\subseteq B $ be finite substructures of $ M. $
	\begin{itemize}
		\item[(i)] If $ B $ is a weak intrinsic extension of $ A $, then $ \delta(B/A)<0. $ Moreover, if $ \intrext{B} $ and $ \delta(B/A)<0 $, then $ B $ contains an increasing chain of weak intrinsic extensions as an initial segment of its constituting chain.
		
		\item[(ii)] If $ B $ is an (weak) intrinsic extension of $ A $, then no proper subset of $ B $ containing $ A $ is (weakly) closed in $ B $. 
		
		\item[(iii)]
		If $ \intrext{B} $ and $ C $ is a finite substructure of $ M $ with $ \intrext{C}, $ then $ \intrext{BC}. $ Moreover, we have that $ \intrext[B]{BC} $ and $ \intrext[C]{BC}. $
		
		\item[(iv)]
		If $ \intrext{B} $ and $ C $ is a finite substructure of $ M $ with $ A\subseteq C, $ then $ \intrext[AC]{BC}. $
		
		\item[(v)] Parts (iii) and (iv) still hold if we replace all occurrences of $ \intrext[]{} $ by $ \wintrext[]{} $.
		
	\end{itemize}
\end{fact}

\begin{dfn}\label{dfnClosure}
	Let $ A\subsetfinite M. $
	\begin{itemize}
		\item[(i)] 
		The \textit{closure} of $ A $ in $ M $ is defined as the union of all intrinsic extensions of $ A $ in $ M $ and is denoted by $ \cl_{M}(A) $.
		
		\item[(ii)] 
		The \textit{weak closure}, or \textit{w-closure} of $ A $ in $ M $, denoted by $ \wcl_{M}(A), $ is the union of all weak intrinsic extensions of $ A $ in $ M $.
		
		\item[(iii)] 
		If $ N\subseteq M, $ the closure of $ N $ in $ M $ is defined as the union of the closures of all finite substructures of $ N, $ namely
		\[ \cl_{M}(N):=\bigcup_{A\subsetfinite N}\cl_{M}(A). \]
		The w-closure of $ N $ in $ M $ is defined in a similar way.
		
		\item[(iv)] The structure  $ N\subseteq M $ is called closed/w-closed in $ M, $ denoted  by $ N<M $ and $ N\leq M $ respectively, if the closure/w-closure of $ N $ in $ M $ is equal to $ N. $
	\end{itemize}
\end{dfn}
\begin{rem}
	Part (ii) of Fact \ref{factIntrExt} shows that if $ B $ is an (weak) intrinsic extension of $ A $, then the (weak) closure of $ A $ in $ B $ equals $ B $. Moreover, the closure/w-closure of a substructure $ N\subseteq M $ is the smallest subset of $ M $ that contains $ N $ and is closed/w-closed in $ M. $
\end{rem}

To obtain the suitable properties required for the machinery of amalgamation classes we restrict ourselves to the following subclass of finite $ \mathcal{L} $-structures.\[ \kpluszero:=\bigg\{A\Big| A \text{ is a finite } \mathcal{L}\text{-structure with } \varnothing<A\bigg\}. \]
Also, $ \kpluszerobar $ denotes the class of $ \mathcal{L} $-structures $ M $ whose finite substructures belong to $ \kpluszero. $

\begin{notation} 
	Let $ A\subsetfinite M\in\kpluszerobar $ and $ A\subseteq B. $
	\begin{itemize}
		\item[(i)] 
		The number of distinct realizations or copies of $ B $ over $ A $ in $ M $ is denoted by $ \chi_{M}(B/A). $ 
		
		\item[(ii)] 
		The maximal number of the disjoint realizations or copies of $ B $ over $ A $ in $ M $ is denoted by $ \chi_{M}^{*}(B/A). $
		
	\end{itemize} 
\end{notation}
The value of $ \chi_{M} $ and $ \chi_{M}^{*} $ is considered to be either a natural number or $ \infty. $

The following fact is somewhat a rephrasing form of Lemma 3.19 in \cite{BaldwinShi-StableGen}:
\begin{fact}\label{factUpperBound}
	For any $ A\in\kpluszero $ and for every structure $ M\in\kpluszerobar $ with $ A\subseteq M, $ if $ B $ is a weak intrinsic extension of $ A $, then we have that
	\[ \chi_{M}(B/A)\leq \frac{\delta(A)}{\card[R^{B/A}]}\leq\delta(A)\leq\card[A],\]
	where $ R^{B/A} $ denotes all the relations with at least one of their components belonging to $ A $ and at least one component belonging to $ B\backslash A. $
	
\end{fact}

Fact \ref{factUpperBound} shows that in any structure $ M\in\kpluszerobar $ the number of realizations of a weak intrinsic extension over the base set is uniformly bounded above; in particular, we have the following fact which is well known in the literature:

\begin{fact}\label{factWClosureFinite}
	Suppose that $ A\subsetfinite M\in\kpluszerobar. $ Then $ \wcl_{M}(A) $ is finite and is a subset of the algebraic closure of $ A $ in $ M. $
\end{fact}

\begin{lem}\label{lmaInfiniteChiImplies}
	If $ A\subsetfinite M\in\kpluszerobar, \intrext{B} $ and $ \chi_{M}^{*}(B/A)=\infty $, then $ \delta(B/A)=0. $ 
\end{lem}
\begin{proof}
	Part (i) of Fact \ref{factIntrExt} implies that if $ \delta(B/A)<0 $, then $ B $ contains a weak intrinsic extension of $ A $, say $ B' $, as an initial segment of its constituting chain over $ A $ (Definition \ref{dfnIntrExt}). By Fact \ref{factWClosureFinite}, $ B' $ is part of the algebraic closure of $ B $ in $ M $ which justifies the claim of the Lemma.
	
\end{proof}

The observation below can easily be seen, but for better referencing we state it as a lemma:
\begin{lem}\label{lmaWclOfCl}
	If $ A\subsetfinite M\in\kpluszerobar $ is weakly closed, then any intrinsic extension of $ A $ in $ M $ is also weakly closed.
\end{lem}
\begin{proof}
	Suppose that $ B $ is an intrinsic extension of $ A $ that is not weakly closed. Then, there exists a finite structure $ C\subsetfinite M $ satisfying $ \delta(C/B)<0. $ Hence, we have that 
	\[ \delta(C/A)=\delta(C/B)+\delta(B/A)<\delta(B/A)\leq 0 \]
	contradicting the assumption that $ A $ is weakly closed in $ M $.
\end{proof}
\begin{lem}\label{lmaChiStarFreeAmalgam}
	Let $ B $ and $ C $ be distinct intrinsic extensions of $ A. $ In any structure $ M\in\kpluszerobar $ containing $ A $ as a substructure with 
	\[ \chi_{M}^{*}(B/A)=\infty\quad\text{and}\quad\chi_{M}^{*}(C/A)=\infty, \]
	there exist infinitely many disjoint copies of $ B $ and $ C $ that are mutually in free amalgamation over $ A. $
\end{lem}
\begin{proof}
	By Lemma \ref{lmaInfiniteChiImplies}, we know that $ \delta(B/A)=\delta(C/A)=0. $ Fix a copy of $ B $ over $ A, $ say $ B_{1} $ and note that for another copy of $ B $ that is disjoint from $ B_{1} $ over $ A, $ say $ B_{2}, $ if there are relations preventing $ B_{2} $ from being freely amalgamated with $ B_{1} $ over $ A, $ then we will have that
	\[ \delta(B_{2}/AB_{1})<\delta(B_{2}/A)=0. \]
	This shows that there exist at most $ \delta(B_{1}) $-many of such copies over $ A $, otherwise, a finite structure with a negative predimension will appear within $ M $ which is impossible due to the fact that $ M\in\kpluszerobar $. Therefore, there are infinitely many disjoint copies of $ B $ that are in free amalgamation with $ B_{1} $ over $ A. $ Iterating this argument leads to finding infinitely many copies of $ B $ that are mutually in free amalgamation over $ A; $ let $ \{B_{i}\}_{i\in\omega} $ enumerate them.
	
	Now, by fixing a copy of $ C, $ say $ C_{1}, $ and considering the fact that the structures $ B_{i} $ are mutually disjoint over $ A $, we obtain an infinite subsequence of $ \{B_{i}\}_{i\in\omega} $ that are disjoint from $ C_{1} $ over $ A $. Again, this argument can be iterated in order to find the desired disjoint copies in $ M. $
\end{proof}
\begin{notation}\label{notKPlusFirstOrder} 
	The first-order theory $ T_{0} $ consists of all formulas of the form $ \forall\bar{x}\neg\diag_{A}(\bar{x}) $ where $ A $ ranges over all finite structures $ A\not\in\kpluszero. $
\end{notation}

\textbf{A historical remark.} By generalizing Fra\"iss\'e amalgamation methods, the machinery of Hrushovski constructions was first introduced in \cite{Hrushovski-NewStrongly} to refute the Zilber's trichotomy conjecture on the available pregeometries inside an $ \aleph_1 $-categorical structure. More details on these constructions can be found in \cite{Wagner-Relational} and \cite{BaldwinShi-StableGen}. The next subsection contains an expansion of the language which has its origins in Hrushovski's results appeared in \cite{Hrushovski-SimpLascarGroups} and also the machinery used in \cite{Pourmahd-SmoothClasses} and \cite{Pourmahd-SimpleGen}.

\vspace*{.3cm}
\subsection{A Natural Expansion of the Language}\label{subsecNat}\hfil


We will define an expansion by definition $ \mathcal{L}^{+} $ of $ \mathcal{L} $ to reinforce the notion of a substructure by making it capable of preserving a notable part of closures in $ \mathcal{L}^{+} $-extensions.
More precisely, given two arbitrary $\mathcal{L}^{+}$-structures $ M\subseteq_{+} N, $ the new augmented predicate symbols together with their supporting axioms would make every finite part of $ \cl_{N}(M) $ realizable in $ M $; in particular, $ M\subseteq_{+} N $ would imply that $ M\leq N $.

To define $ \mathcal{L}^{+}, $ for any pair of structures $ (A,B)\in\kpluszero\times\kpluszero $ with $ B $  an intrinsic extension of $ A $, we add a new predicate $ R_{(A,B)} $ to the language. Namely, let
\[ \mathcal{L}^{+}:=\mathcal{L}\cup\bigg\{R_{(A,B)}(\bar{x})\Big|\intrext{B}\in\kpluszero , \card[\bar{x}]=\card[A]\bigg\}. \]
The $ \mathcal{L}^{+} $-theory $ T_{\nat} $ strengthens $ T_0 $ by adopting the following inductive axioms that makes the new predicates be interpreted as naturally expected:
\[ T_{\nat}:=T_{0}\cup\bigg\{\forall\bar{x}\Big[R_{(A,B)}(\bar{x})\leftrightarrow\big(\diag_{A}(\bar{x})\wedge\exists\bar{y}\diag_{B}(\bar{x};\bar{y})\big)\Big]\bigg|\intrext[A]{B}\in\kpluszero \bigg\}. \]
\begin{rem}
	Any $\mathcal{L}$-structure $ M\in\kpluszerobar $ can be uniquely expanded to an $ \mathcal{L}^{+} $-structure that is a model of $ T_{\nat}. $ But, to ease notation and in dealing with models of $ T_{\nat} $ we identify a model $ M $ with its expansion.
	
\end{rem} 
\begin{lem}\label{lmaLplusLeqstar}
	Suppose that $ M,N\models T_{\nat}. $
	\begin{itemize}
		\item[(i)] If $ M<N, $ then $ M\subseteq_{+} N. $
		
		\item[(ii)] If $ A\subsetfinite M\subseteq_{+} N $ and $ \intrext{B} $ with $ \chi_{N}(B/A)<\infty, $ then we have that
		\[ \chi_{M}(B/A)=\chi_{N}(B/A) \]
		and every copy of $ B $ over $ A $ in $ N $ lies inside $ M. $ Moreover, if $ \chi_{N}(B/A)=\infty, $ then $ \chi_{M}(B/A)=\infty. $
		
		\item[(iii)] If $ M\subseteq_{+} N, $ then $ M\leq N. $
		
	\end{itemize}
\end{lem}
\begin{proof}
	Part (i) directly follows from the definitions. For part (ii), by Fact \ref{factIntrExt}, any finite union of intrinsic extensions is again an intrinsic extension. This, using the hypothesis that $ M\subseteq_{+} N $ implies the result. For part (iii), by Fact \ref{factUpperBound}, for any $ A\subsetfinite M $ and any intrinsic extension $ B $ of $ A $ with $ \delta(B/A)<0 $ we have that
	\[ \chi_{M}(B/A)<\infty. \]
	Hence, by part (ii) of this lemma, we have $ B\subseteq M $ proving $ M\leq N. $
\end{proof}
A more detailed account of the idea of expanding $ \mathcal{L} $ to $ \mathcal{L}^{+} $ can be found in \cite{Pourmahd-SmoothClasses} and \cite{Pourmahd-SimpleGen} where it was used to study Hrushovski constructions in the context of simplicity.

\begin{lem}\label{lmaAmalgam}\hfill
	
	\begin{itemize}
		\item[(i)] 
		Let $ \langle\mathcal{K},\sqsubseteq\rangle\in\big\{\langle\kpluszerobar,\leq\rangle, \langle\kpluszerobar,<\rangle \big\}$. Then $ \langle\mathcal{K},\sqsubseteq\rangle $ has the full amalgamation property; i.e. if $ M_0, M_1 $ ‌and $ M_2 $ are $ \mathcal{L} $-structures in $ \mathcal{K} $ with $ M_0\subseteq M_1 $ and $ M_0\sqsubseteq M_2 $, then $ M=\freejoin{M_1}{M_0}{M_2} $ belongs to $ \mathcal{K} $, and $ M_1\sqsubseteq M. $ Furthermore, if $ M_0\leq_n M_1 $ then $ M_2\leq_n M $.
		
		\item[(ii)] 
		$ \big\langle\hspace*{-2pt}\Mod(T_{\nat}),\subseteq_{+}\hspace*{-4pt}\big\rangle $ has the full amalgamation property.
		
	\end{itemize}
	
\end{lem}
\begin{proof}
	Part (i) is well known in the literature.
	
	For part (ii), Suppose that $ M_{0}\subseteq M_{1},
	M_{0}\subseteq_{+} M_{2} $ and let $ M:=\freejoin{M_{1}}{M_{0}}{M_{2}}. $ Part (iii) of Lemma \ref{lmaLplusLeqstar} implies that $ M_0\leq M_2 $ and hence, by part (i), $ M $ belongs to $ \kpluszerobar. $
	
	To show that $ M_{1}\subseteq_{+} M, $ let $ \intrext{B} $ with $ A\subseteq M_{1} $ and $ B\subseteq M. $ Since $ \intrext[B\cap M_1]{B} $ and $ B\cap M_1 $ is already a subset of $ M_1 $, without loss of generality, we may assume that $ B\cap M_{1}=\varnothing $ and $ B\subseteq M_2\backslash M_1 $. Let $ A_{0}:=A\cap M_{0} $ and $ A_{1}:=A\backslash M_{2}. $ Since $ B $ is in free amalgamation with $ A\backslash M_{2} $ over $ A\cap M_{0}, $ we have that $ \intrext[A_{0}]{B}. $ 
	
	If $ \chi_{M_{2}}(B/A_{0})<\infty, $ part (ii) of Lemma \ref{lmaLplusLeqstar} shows that $ B\subseteq M_{0}. $ 
	
	In the case that  $ \chi_{M_{2}}(B/A_{0})=\infty, $ by part (ii) of  Lemma \ref{lmaLplusLeqstar}, we have that $ \chi_{M_{0}}(B/A_{0})=\infty. $ Since $ B $ is an intrinsic extension of $ A_0, $ the relative predimension of $ B $ over $ A_{0} $ is less than or equal to zero. Hence, for a given copy of $ B $ in $ M_{0} $ that is not in free amalgamation with $ A_{1} $ over $ A_{0}, $ say $ B' $,  we have that
	\[ \delta(B'/A)<\delta(B'/A_{0})\leq 0 .\]
	Therefore, by Fact \ref{factUpperBound}, the number of distinct copies of such a structure $ B' $ over $ A $ in $ M_{0} $ is bounded by $ \delta(A) $. 
	
	On the other hand, considering the diagram of $ B $ over $ A, $ there are only finitely many possibilities available for $ B $ being not in free amalgamation with $ A_{1} $ over $ A_{0}. $ Therefore, there exist infinitely many copies of $ B $ inside $ M_{0} $ that are in free amalgamation with $ A_{1} $ over $ A_{0}. $ This, in particular, provides a copy of $ B $ over $ A $ in $ M_{1} $ which shows that $ M_{1}\subseteq_{+}M. $
\end{proof}

\vspace*{.3cm}
\section{Introducing $\tmu$ }\label{secTmu}
By fixing a certain natural-valued function $ \mu $ defined on pairs of intrinsic extensions, we consider a class of infinite structures called $ \mu $-\textit{paracollapsed}, or more briefly \textit{paracollapsed} structures. This class of structures will be axiomatized by an $\mathcal{L}$-theory $ \tmuinfinity $. Thereafter, we introduce an $ \mathcal{L}^{+} $-theory $ \tmu $ that will eventually axiomatize the class of existentially closed models of $ \tmuinfinity\cup T_{\nat} $. In subsection \ref{subsecCons}, we prove that $ \tmu $ is consistent by showing that every existentially closed model of $ \tmuinfinity\cup T_{\nat} $ is a model of $ \tmu $. In subsection \ref{subsecQE} we show that $ \tmu $ eliminates quantifiers in $ \mathcal{L}^{+} $. 

Axioms of $ \tmu $ are obtained by a suitable adaptation of what is called \textit{semigenericity} in the literature (see \cite{BaldwinShelah-Randomness}) together with strengthening it by a new axiom scheme that will be called the axiom of \textit{intrinsic semigenericity}.

\begin{convention}
	From now on, any finite structure under question would belong to $ \kpluszero. $
\end{convention}
\begin{assumption}\label{assummMu}
	We fix a function $ \mu:\omega\backslash\{0\}\times\omega\backslash\{0\}\longrightarrow\omega $ with the following properties:
	\begin{itemize}
		\item[(i)] 
		For every $ m, n\in\omega\backslash\{0\} $ we have that $ \mu(m,n)\geq m. $

		\item[(ii)] 
		For any nonzero natural numbers $ m_1,m_2,n_1 $ and $ n_2 $ we have that 
		\[ \mu(m_1+m_2,n_1+n_2)>\mu(m_1,n_1)+\mu(m_2,n_2). \]
		

		%
		%
	\end{itemize}
\end{assumption}
\begin{rem}
	Part (ii) of assumption above ensures that for any $ m_1, n_1, m_2 $ and $ n_2 $ in the domain of $ \mu $, whenever $ m_1\leq m_2 $ and $ n_1\leq n_2 $ we have that $ \mu(m_1,n_1)\leq\mu(m_2,n_2) $. An example of $ \mu $ can be obtained using the real-valued function that maps each pair of nonzero natural numbers $ (m,n) $ to $ e^{m+n} $. 
\end{rem}
\begin{dfn}\label{dfnUCC}
	A structure $ M\in\kpluszerobar $ is said to be $ \mu $-\textit{paracollapsed}, or briefly \textit{paracollapsed}, if for every pair of finite structures $ (A,B) $ with $ A\subsetfinite M $ and $ \intrext{B} $ either we have 
	\[ \chi^{*}_{M}(B/A)\leq\mu(\card[A],\card[B]), \]
	or $ \chi^{*}_{M}(B/A)=\infty. $ The class of all paracollapsed structures is denoted by $ \kmuinfinity. $
\end{dfn}

\begin{rem}\label{rmkUCCFirstOrder}
	The following first order extension of $ T_{0} $ axiomatizes the class $ \kmuinfinity: $
	\[ \tmuinfinity:=T_{0}\cup\bigg\{\forall\bar{x}\Big[\diag_{A}(\bar{x})\wedge\exists^{\geq n}\bar{y}\diag_{B}(\bar{x};\bar{y})\rightarrow\exists^{\geq n+1}\bar{y}\diag_{B}(\bar{x};\bar{y})\Big]\bigg\}, \]
	where $ (A,B) $ ranges over all pairs of finite structures with $ \intrext{B} $, and $ n $ ranges over all natural numbers strictly greater than $ \mu(\card[A],\card[B]). $
\end{rem}
\begin{notation}
	We denote the $ \mathcal{L}^{+} $-theory $ \tmuinfinity\cup T_{\nat} $ by $ \tplus. $
\end{notation}
\begin{dfn}\label{dfnOmittable}
	Let $ B\subseteq C $ be finite.
	\begin{itemize}
		\item[(i)] 
		A model $ N\models\tmuinfinity $ containing an isomorphic copy of $ B, $ say $ \bar{b}, $ is said to \textit{omit} $ C $ over $ \bar{b} $ if 
		\[ N\models\neg\exists\bar{z}\diag_{C}(\bar{b};\bar{z}). \]
		
		\item[(ii)]
		$ C $ is called $ \mu $-\textit{paracollapse-omittable} over $ B, $ or briefly \textit{omittable}, if there exists a model $ N\models\tmuinfinity $ with an embedding $ f:B\longrightarrow N $ that omits $ C $ over $ fB. $
		
	\end{itemize}
	
\end{dfn}

Before proceeding further, in the remark below we introduce a set of first order formulas that will be used in the rest of this section:
\begin{rem}\label{rmkFormulas}\hfil
	\begin{itemize}
		\item[(i)]
		Given a finite structure $ A $ and a natural number $ n, $ there is a formula $ \sigma_{\scalebox{.5}{$ (A,n) $}}(\bar{x}) $ with $ \card[\bar{x}]=\card[A] $ such that for any $ M\in\kpluszerobar $ and any tuple $ \bar{a}\in M $ we have that  $ M\models\sigma_{\scalebox{.5}{$ (A,n) $}}(\bar{a}) $ if and only if $ \bar{a}\cong A $ and $ \bar{a}\leq^{n}M $. Namely, $ \sigma_{\scalebox{.5}{$ (A,n) $}}(\bar{x}) $ is the following formula
		\[ \diag_{A}(\bar{x})\wedge\bigwedge_{B\in\mathcal{K}_{\scalebox{.5}{$ (A,n) $}}}\neg\exists\bar{y}_{_{\scalebox{.5}{B}}}\diag_{B}(\bar{x};\bar{y}_{_{\scalebox{.5}{B}}}), \]
		where $ \mathcal{K}_{\scalebox{.5}{$ (A,n) $}} $ is the collection of all finite structures $ B\in\kpluszero $ with $ \card[B\backslash A]\leq n $ and $ A\not\leq B. $
		
		\item[(ii)]
		Let $ B $ be a finite structure and $ \Omega $ a finite set of omittable structures over $ B. $ Then, there is a first-order formula, as below, denoted by $ \omega_{\scalebox{.5}{$ (B,\Omega) $}}(\bar{y}) $ with $ \card[\bar{y}]=\card[B] $ which claims in any model $ N\models\tmuinfinity $ with $ N\models\omega_{\scalebox{.5}{$ (B,\Omega) $}}(\bar{b}) $ that the tuple $ \bar{b} $ is isomorphic to $ B $ and that $ N $ omits each $ C\in\Omega $ over $ \bar{b}: $
		\[ \diag_{B}(\bar{y})\wedge\bigwedge_{C\in\Omega}\neg\exists\bar{z}_{\scalebox{.4}{$ C $}}\diag_{C}(\bar{y};\bar{z}{\scalebox{.4}{$ C $}}). \]
	\end{itemize}
	
\end{rem}
\begin{dfn}\label{dfnUCCGeneric}
	Let $ A\subsetfinite M\in\kmuinfinity $ and $ n\in\omega $. Also, suppose that $ B\in\kpluszero $ with $ A\subseteq B $ and $ \Omega $ is a set of omittable structures over $ B $ such that each $ C\in\Omega $ is an intrinsic extension of $ B $ while $ A(C\backslash B) $ is not an intrinsic extension of $ A. $
	\begin{itemize}
		\item[(i)]
		$ M $ is said to satisfy the axiom of \textit{semigenericity} for $ A, B, \Omega $ and $ n $, if whenever $ A\leq_{n} M $ and $ A<B $ there exists an embedding $ f:B\longrightarrow M $ with $ f\restriction_{A}=\identity_{A} $ and $ fB\leq_{n}M $ such that $ M $ omits each $ C\in\Omega $ over $ fB. $ That is, $ M $ satisfies the following formula
		\[\hspace*{-75pt}\opaxsem_{(\scalebox{.55}{$ A,B,\Omega,n $})}:\qquad\qquad \forall\bar{x}\bigg[\sigma_{\scalebox{.5}{$ (A,n) $}}(\bar{x})\rightarrow\exists\bar{y}\Big(\sigma_{\scalebox{.5}{$ (B,n) $} }(\bar{x};\bar{y})\wedge\omega_{\scalebox{.5}{$ (B,\Omega) $}}(\bar{x}\bar{y})\Big)\bigg]. \]
		
		\item[(ii)]
		We say that $ M $ satisfies the axiom of \textit{intrinsic semigenericity} for $ A, B, \Omega $ and $ n $, if whenever $ A\leq_{n} M $ and $ \intrext[A]{B} $ with $ \chi^{*}_{M}(B/A)=\infty $, there exists an embedding $ f:B\longrightarrow M $ with $ f\restriction_{A}=\identity_{A} $ and $ fB\leq_{n}M $ such that $ M $ omits each $ C\in\Omega $ over $ fB. $ In other words, $ M $ satisfies the following formula
		\begin{align*}
			\hspace*{-25pt}\opaxisem_{(\scalebox{.55}{$ A,B,\Omega,n $})}:\hspace*{7pt} \forall\bar{x}\bigg[\Big(\sigma_{\scalebox{.5}{$ (A,n) $}}(\bar{x})\wedge \exists^{\raisebox{3pt}{\scalebox{.4}{$ \infty $}}}\bar{y}\chi_{\scalebox{.5}{$ (A,B) $}}^{*}(\bar{x};\bar{y})\Big)\rightarrow\exists\bar{y}\Big(\sigma_{\scalebox{.5}{$ (B,n) $}}(\bar{x};\bar{y})\wedge\omega_{\scalebox{.5}{$ (B,\Omega) $}}(\bar{x}\bar{y})\Big)\bigg].
		\end{align*}
	\end{itemize}
\end{dfn}

\begin{dfn}\label{defUCCTheory}
	$ \tmu $ denotes the following $ \Pi_{2} $-axiomatizable $ \mathcal{L}^{+} $-theory 
	\[ \tplus\cup\Big\{\opaxsem_{(\scalebox{.55}{$ A,B,\Omega,n $})}\Big\}\cup\Big\{\opaxisem_{(\scalebox{.55}{$ A,B,\Omega,n $})}\Big\}, \]
	where $ A, B, \Omega $ and $ n $ are correspondingly as mentioned in Definition \ref{dfnUCCGeneric}.
\end{dfn}

\vspace*{.3cm}
\subsection{Consistency of $\tmu$}\label{subsecCons}\hfill

We show that any existentially closed model of $ \tplus $ is a model of $ \tmu $. Towards this goal, Lemma \ref{lmaUCCExtension} forms a major step by enabling us to extend canonically a given structure $ M\in\kpluszerobar $ to a model of $ \tmuinfinity $, called the $ \mu $-\textit{paracollapsed hull} of $ M $. As a stronger consequence, in Proposition \ref{propModelCompletion} below, we will see that $ \tmu $ is actually the model completion of the theory $ \tplus $.

\begin{lem}\label{lmaUCCExtensionFirstStep}
	Given a countable $ M\models T_{0}, $ there exists a countable model $ \overline{M}\models T_{0} $ extending $ M $ with $ M\leq \overline{M} $ and having the following properties:
	\begin{itemize}
		\item[(i)] 
		For any $ \intrext{B} $ with $ A\subsetfinite M $ and $ \chi_{M}^{*}(B/A)>\mu(\card[A],\card[B]) $ we have that 
		\[ \chi_{\overline{M}}^{*}(B/A)=\infty. \]
		
		\item[(ii)] 
		For any $ N\models\tmuinfinity $ and any embedding $ f:M\longrightarrow N $ with $ fM\leq N, $ there exists an embedding $ \bar{f}:\overline{M}\longrightarrow N $ extending $ f $ such that $ \bar{f}\overline{M}\leq N. $
		
		\item[(iii)]
		For any $ M'\subseteq M $ with $ M'\models\tmuinfinity, $ if $ M'\subseteq_{+} M $ ($ M'<M $) then we have that $ M'\subseteq_{+} \overline{M} $ ($ M'<\overline{M} $).
	\end{itemize}
\end{lem}
\begin{proof}
	Let $ \big\{(A_{n},B_{n}):n\geq 1\big\} $ enumerate all pairs that violate the paracollapsed condition in $ M; $ this means that for each $ n\geq 1 $ the structure $ B_{n} $ is an intrinsic extension of $ A_{n}, $ and $ \chi_{M}^{*}(B_{n}/A_{n}) $ is finite and is strictly greater than $ \mu(\card[A_{n}],\card[B_{n}]). $ 
	
	Let $ M_{0}:= M $ and for $ n\geq 1 $ define $ M_{n} $ as the free amalgamation of $ M_{n-1} $ (over $ A_{n} $) with infinitely many disjoint copies of $ B_{n} $ that are mutually in free amalgamation over $ A_{n}. $ Finally, define $ \overline{M} $ as the union of all structures $ M_{n} $. It is obvious that each $ M_{n} $ and hence $ \overline{M} $ belongs to $ \kpluszerobar. $ Also, condition (i) in the lemma is obvious by our way of constructing $ \overline{M} $ and condition (ii) follows from the construction and Lemma \ref{lmaChiStarFreeAmalgam}.
	
	For condition (iii), suppose that $ M'\subseteq_{+} M. $ For every $ n\geq 1 $ we show that $ M'\subseteq_{+} M_{n}. $
	
	Let $ \bar{A}_{n}:=A_{n}\cap M' $ and $ L_{n}:=\cl_{B_{n}}(\bar{A}_{n}). $ The case of $ L_{n}=\bar{A}_{n} $ leads trivially to the stronger property of $ M'< M_{n}. $ Hence, we suppose that $ L_{n}\backslash\bar{A}_{n}\neq\varnothing. $
	
	Since $ \chi_{M}^{*}(B_{n}/A_{n}) >\mu(\card[A_{n}],\card[B_{n}]), $ we have that 
	\[  \chi_{M}^{*}(L_{n}/\bar{A}_{n})\geq \chi_{M}^{*}(B_{n}/A_{n}) >\mu(\card[A_{n}],\card[B_{n}]).  \]
	This holds because any disjoint copy of $ B_{n} $ over $ A_{n} $ carries within itself a disjoint copy of $ L_{n} $ over $ \bar{A}_{n} $. 
	
	On the other hand, by part (ii) in Assumption \ref{assummMu}, we know that $ \mu $ is an increasing function according to the lexicographic order. In particular, we have that $ \mu(\card[A_{n}],\card[B_{n}])\geq \mu(\card[\bar{A}_{n}],\card[L_{n}]) $ which by the above inequality implies that
	\[ \chi_{M}^{*}(L_{n}/\bar{A}_{n})>\mu(\card[\bar{A}_{n}],\card[L_{n}]). \]
	Since $ \bar{A}_{n}\subseteq M' $ and $ M'\subseteq_{+} M, $ there exist at least $ \Big(\mu(\card[\bar{A}_{n}],\card[L_{n}])+1\Big) $-many disjoint copies of $ L_{n} $ over $ \bar{A}_n $ inside $ M'. $ This, using the fact that $ M'\models\tmuinfinity, $ leads to finding infinitely many disjoint copies of $ L_{n} $ over $ \bar{A}_{n} $ in $ M'. $ This shows that in the process of constructing $ M_n $ we did not add any new finite intrinsic extension over $ M' $ that had been already absent in $ M' $; this proves that $ M'\subseteq_{+} M_{n}. $
	
	The hypothesis of $ M'<M $ is reduced to the case of $ L_{n}=\bar{A}_{n}. $
\end{proof}
The following corollary will be used in the proof of Theorem \ref{thmQE}:
\begin{crl}\label{crl-ommitable}
	For two finite structures $B\subseteq C\in\kpluszero$, there exists an efficient procedure to detect if $C$ is omittable over $B$.
\end{crl}
\begin{proof}
	Proceeding as in the proof of \Cref{lmaUCCExtensionFirstStep}, it is easy to verify that passing the following procedure for $(B,C)$ is equivalent to $C$ being omittable over $B$:
	
	\begin{itemize}
		\item[(1)] Enumerate all pairs $(A,D)$ of substructures of $B$ with $\intrext{D}$ and such that 
		\[ \chi^{*}_{B}(D/A)>\mu(\card[A],\card[D]). \]
		
		\item[(2)] Divide $C$ into finitely many components $C_1,\ldots, C_k$ so that they are all mutually in free amalgamation over $B$. Note that $k$ might be equal to $1$.
		
		\item[(3)] For each $C_i$, let $n=\card[C_i]$ and enumerate all sequence of structures $E_1,\ldots, E_n$ where each $E_j$ is isomorphic, over $B$, to one of the structures $D$ obtained in the first step. Note that some of the structures $E_j$ might be isomorphic. For each such sequence, freely amalgamate all $E_j$ over $B$. 
		
		\item[(4)] For each $C_i$, check if it appears as a substructure of the structures built in the previous step. If not, then $C_i$, and hence $C$, is omittable over $B$. 
		
		\item[(5)] If for each $i=1,\ldots, k$, the answer to the previous step is not negative, then $C$ is not omittable over $B$.
	\end{itemize}
	
\end{proof}
\begin{lem}\label{lmaUCCExtension} Given a countable $ M\models T_{0}, $ there exists a countable extension $ \mathbbm{M}\models\tmuinfinity $ with $ M\leq\mathbbm{M} $ that satisfies the properties (ii) and (iii) mentioned in Lemma \ref{lmaUCCExtensionFirstStep}. In particular, for any $ \bar{b}\in M $ and $ C\in\kpluszero $ that is omittable over $ \bar{b} $ we have that
	\[ \mathbbm{M}\models\neg\exists\bar{z}\diag_{C}(\bar{b};\bar{z}). \]
\end{lem}
\begin{proof}
	Let $ M_{0}:=M $ and for each $ n\geq 1 $ define $ M_{n} $ to be the structure $ \overline{M}_{n-1} $ obtained by applying Lemma \ref{lmaUCCExtensionFirstStep} for $ M_{n-1} $. Finally, let $ \mathbbm{M}:=\bigcup_{n\in\omega}M_{n}. $
	
	Property (i) in Lemma \ref{lmaUCCExtensionFirstStep} guarantees that $ \mathbbm{M} $ belongs to $ \kmuinfinity. $ It is easy to check that property (iii) is preserved under taking a union.
	
	Given a model $ N\models\tmuinfinity $ with an embedding $ f:M\overset{\leq}{\longrightarrow}N,  $ by property (ii) of Lemma \ref{lmaUCCExtensionFirstStep} we can appropriately embed each $ M_{n} $ into $ N $ over $ M_{n-1} $ which ultimately gives the desired embedding of $ \mathbbm{M} $ into $ N. $ The fact that an omittable structure over a finite substructure of $ M $ is not  realized in $ \mathbbm{M} $ is a direct consequence of Definition \ref{dfnOmittable} and property (ii) of \Cref{lmaUCCExtensionFirstStep}.
\end{proof}

\begin{notation}\label{notUCCHull}
	For a countable model $ M\models T_{0}, $ the structure $ \mathbbm{M}\models\tmuinfinity $ obtained in Lemma \ref{lmaUCCExtension} is called the $ \mu $-\textit{paracollapsed hull}, or briefly the \textit{paracollapsed hull}, of $ M $ and is denoted by $ \uccext. $
\end{notation}

\begin{lem}\label{lmaExCloedIsUCCGeneric}
	Any existentially closed model of $ \tplus $ is a model of $ \tmu $.
\end{lem}
\begin{proof}
	Let $ M $ be a countable existentially closed model of  $ \tplus. $ To show that $ M $ satisfies the axiom \axsem\hspace*{-3pt}, suppose that $ \bar{a}<\bar{a}\bar{b}\in\kpluszero ,\bar{a}\leq_{n} M, $ and let $ \Omega $ be a finite set of omittable structures over $\bar{a}\bar{b}.$
	
	If $ N_0 $ denotes the free amalgamation of $ M $ with $ \bar{a}\bar{b} $ over $ \bar{a} $, by Lemma \ref{lmaAmalgam}, we know that $ N_{0} $ belongs to $ \kpluszerobar $ and moreover $ M<N_{0}. $ Let $ N $ be the paracollapsed hall of $ N_0 $, namely $ N=\uccext[N_{0}]$.
	
	By Lemma \ref{lmaUCCExtension}, for each $ C\in\Omega $ we have that
	\[ N\models\neg\exists\bar{z}\diag_{C}(\bar{a}\bar{b};\bar{z}). \]
	By Remark \ref{rmkFormulas}, this actually means that
	\[ N\models \omega_{\scalebox{.5}{$ (\bar{a}\bar{b},\Omega) $}}(\bar{a}\bar{b}). \]
	
	On the other hand, by Lemma \ref{lmaAmalgam}, we have that $ \bar{a}\bar{b}\leq_{n} N_{0} $ which, using Lemma \ref{lmaUCCExtension}, implies that $ \bar{a}\bar{b}\leq_{n} N. $ This, in particular, means that for any $ n $ we have
	\[ N\models\sigma_{\scalebox{.5}{$ (\bar{a}\bar{b},n) $} }(\bar{a}\bar{b}). \]
	
	By property (iii) in Lemma \ref{lmaUCCExtension}, we have that $ M $ is closed in $ N $. This, by part (i) of Lemma \ref{lmaLplusLeqstar}, implies that $ M\subseteq_{+} N. $ 
	
	By Lemma \ref{lmaUCCExtension}, we know that $ N $ is a model of $ \tmuinfinity $. This model can readily be expanded to a model of $ T_{\nat} $, and hence a model of $ \tplus $. Moreover, by the hypothesis, $ M $ is an existentially closed model of $ \tplus $ which altogether yields
	\[ M\models\opaxsem_{(\scalebox{.55}{$ \bar{a},\bar{a}\bar{b},\Omega,n $})}. \]
	
	The proof of $ M\models\opaxisem_{(\scalebox{.55}{$ \bar{a},\bar{a}\bar{b},\Omega,n $})} $ proceeds using a similar argument; the only difference is that the condition of $ M\subseteq_{+}N_{0} $ follows directly from the hypothesis of $ \chi_{M}^{*}(\bar{b}/\bar{a})=\infty $ and not from Lemma \ref{lmaAmalgam}.
	
\end{proof}

\begin{thm}\label{thmCons}
	$ \tmu $ is consistent.
\end{thm}
\begin{proof}
	Being an inductive theory, any model of $ \tplus $ can be extended to an existentially closed model. Therefore, by Lemma \ref{lmaExCloedIsUCCGeneric}, any such model is a model of $ \tmu $.    
\end{proof}

The following corollary will be used in the proof of Theorem \ref{thmQE}:
\begin{crl}\label{crlJEP}
	$ \tmu $ has the joint embedding property.
\end{crl}
\begin{proof}
	Suppose that $M_1$ and $M_2$ are two model of $\tmu$. Let $M$ denote the disjoint union of $M_1$ and $M_2$ as a structure in $\kpluszerobar$, then apply \Cref{lmaUCCExtensionFirstStep} and \Cref{lmaUCCExtension} to obtain the paracollapsed hull of $M$. By definition, $\uccext$ is a model of $\tmuinfinity$, which can easily be expanded to a model of $\tplus $, which, in turn, can be extended to an existentially closed model of the same theory (see the proof of \Cref{thmCons}). By \Cref{lmaExCloedIsUCCGeneric}, the latter model, say $N$, is a model of $\tmu$. Condition (iii) addressed in \Cref{lmaUCCExtension} guarantees that $M_1\subseteq_{+} N$ and $M_2\subseteq_{+} N$.
\end{proof}

\begin{rem}\label{rmk-Hrushovski-amalgamation}
	Despite all the subtleties it might involve, it seems possible to build a model of $ \tmu $ completely out of scratch by adapting the usual techniques of Fra\"iss\'e-Hrushovski method into the context of infinite objects introduced here, namely the paracollapsed structures. This means primarily to define a suitable notion of \textit{richness} and then to constructively build a rich model using amalgamation over the structures in $ \kmuinfinity $. In that way, the obtained rich model, supposed to be universal and \textit{ultra-homogeneous} with regard to the structures in $ \kmuinfinity $, would play the role of a model of $ \tmu $. However, as mentioned earlier, we chose a different direction by following Hrushovski's approach in \cite{Hrushovski-SimpLascarGroups} that is based on techniques of Robinson model theory in axiomatizing existentially closed models. 
\end{rem}

\vspace*{.3cm}
\subsection{Completeness and Quantifier Elimination}\label{subsecQE}\hfill

In this subsection we prove that $ \tmu $ is complete and eliminates quantifiers in $ \mathcal{L}^{+}. $ This, using Lemma \ref{lmaExCloedIsUCCGeneric}, will imply that $ \tmu $ is a model companion for $ \tplus $. In fact, we will see in Proposition \ref{propModelCompletion} that $ \tmu $ is actually the model completion of $ \tplus $.

\begin{thm}\label{thmQE}
	$ \tmu $ is complete decidable theory that eliminates quantifiers in $ \mathcal{L}^{+} $.
\end{thm}
\begin{proof}
	To prove quantifier elimination, let $ M_1 $ and $ M_2 $ be two models of $ \tmu $ with a common $ \mathcal{L}^{+} $-substructure $ A=M_1\cap M_2. $ Being a model of $ (T_{\nat})_{\forall} $, the structure $ A $ satisfies the following weak version of axiom schemes in $ T_{\nat} $ for all pair of structures (B,C) with $ \intrext[B]{C} $:
	\[ \forall\bar{x}\Big[\big(\diag_{B}(\bar{x})\wedge\exists\bar{y}\diag_{C}(\bar{x};\bar{y})\big)\rightarrow R_{(B,C)}(\bar{x})\Big]. \]
	
	Suppose that $ \bar{a}\in A $ and $ b $ is an element in $ M_1\backslash A $ satisfying a quantifier free $ \mathcal{L}^+ $-formula $ \varphi(\bar{x}y) $ as
	\[ \bigwedge_{i=1}^{k_1}R_{(\bar{a}b,C_i)}(\bar{x}y) \wedge \bigwedge_{j=1}^{k_2}\neg R_{(\bar{a}b,D_i)}(\bar{x}y),\]
	where for each $ i $ and $ j $ we have that $ \intrext[\bar{a}b]{C_i} $ and $ \intrext[\bar{a}b]{D_j} $. We will show that there exists an element $ b'\in M_2 $ satisfying the same formula. 
	
	By Fact \ref{factWClosureFinite}, the weak closure of a finite subset of $ A $, say $ \bar{c} $, is finite in both $ M_1 $ and $ M_2 $. Using a suitable set of predicates $ R_{(\bar{c},E)}\in\mathcal{L}^{+} $, we can easily show that the weak closure of $ \bar{c} $ in $ M_1 $ is isomorphic to its weak closure in $ M_2 $. Therefore, without loss of generality, we may assume that $ A $ is weakly closed in $ M_1 $ and $ M_2 $; in particular, we may assume that $ \bar{a} $ is weakly closed in both $ M_1 $ and $ M_2 $.
	
	Let
	\[ C=\bigcup_{i=1}^{k_1}C_i. \]
	
	We partition the structures $ D_j $ into two types of structures; recall that any structure $ D_j $ is an intrinsic extension of $ \bar{a}b $. The first type, or \textbf{type-1}, consists of all structures $ D_j $ that are also an intrinsic extension of $ \bar{a} $. All other structures $ D_j $ would belong to the second type of structures, or \textbf{type-2}. 
	
	\vspace*{7pt}
	\textbf{Discussion.}
	If a type-1 structure $ D_j $ is not realized in $ M_1 $ over $ \bar{a} $ at all, then, by the fact that $ M_1 $ is a model of $ T_{\nat}, $ this assumption implies that $ M_1 $ satisfies $ \neg R_{(\bar{a},D_j)}(\bar{a}) $ which easily leads to
	\[ M_2\models \neg R_{(\bar{a},D_j)}(\bar{a}). \]
	Consequently, there would not exist any copy of $ D_j $ over $\bar{a}$ in $ M_2 $. Hence, this kind of type-1 structures are automatically omitted in $ M_2 $ and hence, in the rest of the proof we assume any type-1 structure to be realized in $ M_1 $.
	\hfill \scalebox{.7}{$ \blacksquare_{Discussion} $}
	
	\vspace*{7pt}
	The very formalism of the axioms \axsem and \axisem would guarantee the omission of all type-2 structures over an embedding obtained by those axioms. Therefore, in all the cases below, except for Case B.2, solely the type-1 structures are concerned.
	
	\vspace*{7pt}
	\textbf{Case A:} $ b\in\cl^{*}_{M_1}(\bar{a}). $ In this case, for all the structures $ C $ and $ D_j $ appeared in $ \varphi(\bar{x}y) $ we  have that
	\[ \intrext[\bar{a}]{bC}\quad\text{ and }\quad \intrext[\bar{a}]{bD_j}. \]
	
	If $ bC $ is algebraic over $ \bar{a} $ with $ m $-many conjugates in $ M_1 $, using a suitable set of $ \mathcal{L}^{+} $-predicates, we can easily show the existence of a set $ b'C' $ in $ \acl_{M_2}(\bar{a}) $ that is isomorphic to $ bC $ over $ \bar{a} $ with exactly $ m $-many number of copies in $ M_2 $. Moreover, the existence/non-existence of the structures $ D_j $ over these $ m $-many conjugates of $ bC $ can be described by a finite number of $ \mathcal{L}^{+} $-predicates. Therefore, at least one of the conjugates of $ b'C' $ in $ M_2 $ will satisfy the formula $ \varphi. $
	
	Now suppose that $ bC $ is not algebraic over $ \bar{a} $ in $ M_1 $. 
	To ease notation, we drop the subscript $ j $ from any structure $ D_j $ under consideration below.
	
	For a type-1 structure $ D $, by the Discussion above, we can assume that there is a realization of this structure over $ \bar{a} $ in $ M_1 $; let $ D' $ denote this realization. We prove the following claim for $D$:
	
	\vspace*{7pt}
	\textbf{Claim A.} $ \delta(D/\bar{a}bC)<0. $
	
	\textit{proof of the claim.} First, we show that $ D' $ is in free amalgamation with $ bC $ over $ \bar{a} $. If it was not the case and there were relations violating this free amalgamation, then we would have that 
	\[ \delta(D'/\bar{a}bC)<\delta(D'/\bar{a})\leq 0 \]
	which shows that $ \bar{a}bC $ is not weakly closed in $ M_1 $. By Lemma \ref{lmaWclOfCl}, this would contradict the assumption of $ A\leq M_1 $.
	
	Now, note that the diagram of $ D $ is omitted over $ \bar{a}b $ in $ M_1 $. On the other hand, we just proved that there exists a realization of the diagram of $ D $, over the tuple $\bar{a}$, that is in free amalgamation with $ bC $ over $ \bar{a} $, namely $ D' $. This shows that the diagram of $ D $, over the structure $ \bar{a}bC $, contains at least one relation whose components intersects $ bC $. Hence, we have that
	\[ \delta(D/\bar{a}bC)<\delta(D/\bar{a})=0, \]
	where the latter equality holds because $ D $ is a type-1 structure; that is $ D $ is an intrinsic extension of $ \bar{a} $.
	
	\hfill \scalebox{.7}{$ \blacksquare_{Claim\hspace*{2pt}A} $}
	\vspace*{8pt}
	
	A consequence of the above claim is that any type-1 structure $ D $ with at least one realization over $ \bar{a} $ in $ M_1 $ will be omitted over any sufficiently weakly closed embedding of $ bC $ into $ M_2 $.
	
	Finally, for Case A, let $ n $ denote the maximum of the cardinals  $ \card[D\backslash\bar{a}] $ where $ D $ ranges over all type-1 structures apeearing in $\varphi$. Also, let $ \Omega $ consist of all structures $ D $ of the second type. Now, an application of \axisem for $ \bar{a}, bC, n $ and $ \Omega $ leads to finding an embedding of $ bC $ into $ M_2 $ that omits all type-2 structures $ D $ and that is $ n $-weakly closed in $ M_2 $. The latter property, by the discussion after Claim A, implies that all type-1 structures are also omitted over this embedding.

	\vspace*{7pt}
	\textbf{Case B:} $ b\not\in\cl^{*}_{M_1}(\bar{a})$. Let
	\[ E_0:=\cl^{*}_{bC}(\bar{a}), \quad\text{and}\quad E_1:=bC\backslash E_0. \]
	By the definition of closure, we have that $ \intrext[\bar{a}]{E_0} $ and $ E_0<E_0E_1. $ Our aim would be first to find a suitable copy of $ E_0 $ over $ \bar{a} $ in $ M_2 $, say $ E'_0 $, and then to apply the axiom \axsem in order to embed $ E_1 $ over $ E'_0 $ while omitting all structures $ D_j $ over that embedding.

	As in Case A, we drop the subscripts of the structures $ D_j $ under consideration.
	
	\vspace*{7pt}
	\textbf{Case B.1: $ E_0=\bar{a}. $} By Discussion above, given a type-1 structure $ D $, which is an intrinsic extension of $ \bar{a} $, we can assume that there exists a copy of $ D $ over $\bar{a}$ in $ M_1 $; let $ D' $ denote this realization. We prove the following claim for $D$:
	
	\vspace*{7pt}
	\textbf{Claim B.} $ \delta(D/b)<0. $
	
	\textit{proof of the claim.} First, we show that $ D' $ is in free amalgamation with $ b $ over $ \bar{a} $. If it was not the case and there were relations violating this free amalgamation, $ b $ would belong to the closure of $ \bar{a}D' $ which, using the fact that $ D' $ is already in the closure of $ \bar{a} $, implies that $ b $ belongs to the closure of $ \bar{a} $. This contradicts our assumption of $ b\not\in\cl^{*}_{M_1}(\bar{a}) $. The rest of the proof proceeds exactly similar to the proof of Claim A.
	
	\hfill \scalebox{.7}{$ \blacksquare_{Claim\hspace*{2pt}B} $}
	\vspace*{8pt}
	
	%
	
 In this case, let $ n $ and $ \Omega $ be exactly as they were considered in Case A. Then, an application of \axsem for $ \bar{a}, E_1, n $ and $ \Omega $ leads to finding an embedding of $ E_1 $ into $ M_2 $ that omits all type-2 structures $ D $, and moreover is $ n $-weakly closed in $ M_2 $. The latter property, by Claim B, implies that all type-1 structures are also omitted over this embedding.
	
	\vspace*{7pt}
	\textbf{Case B.2: $ E_0\supsetneq\bar{a} $.} Since $ E_0$ is an intrinsic extension of $ \bar{a} $, the structure $ M_1 $ satisfies $ R_{(\bar{a},E_0)}(\bar{a}) $. Hence, it is easily seen that there exists a copy of $ E_0 $ over $ \bar{a} $ in $ M_2 $; let $ E'_0 $ denote this copy. Lemma \ref{lmaWclOfCl} implies that $ \bar{a}E'_0 $ is weakly closed in $ M_2 $.
	
	We handle all type-1 structures exactly in the same way as we did in Case B.1. For a type-2 structure $ D $, recall that such a structure is an intrinsic extension of $\bar{a}b$ but not an intrinsic extension of $ \bar{a} $. Also, notice that for all structures $ D $, the formula $ \varphi(\bar{x}y) $ only contains information concerning the diagram of $ D $ over $ \bar{x}y $ while remaining silent on the possible interactions between $ D $ and $ C $. In particular, $ \varphi(\bar{x}y) $ contains no information about the diagram of $ D $ over $ E_0 $.
	
	Based on the latter observation, we define a new diagram $ D' $ that extends $ D $ to a new structure over $\bar{a}bE_0$ without adding any new relation; that is, in $ D' $ the structure $ D $ is in free amalgamation with $ E_0 $ over $ \bar{a}b $.
	
	Defined thus, $ D' $ is actually an intrinsic extension of $ \bar{a}bE_0 $ while being not an intrinsic extension over $ \bar{a}E_0 $. Hence, if we let $ n=\card[D\backslash\bar{a}] $ and $ \Omega=\{D'\} $ and then apply \axsem for $ E'_0, E_1, n$ and $ \Omega $, we can find an embedding $ f $ of $ E_1 $ over $ E'_0 $ into $ M_2 $ while omitting $ D' $ over this embedding.
	
	Therefore, the only remaining chance for $ D $ to be realized over $ \bar{a}E'_0f(E_1) $ is to be not isomorphic to $ D' $ by having some extra relations over $ E'_0 $. If $ D'' $ is such a realization, then we have that 
	\[ \delta(D''/\bar{a}E'_0f(E_1))<\delta(D''/\bar{a}f(E_1))\leq\delta(D''/\bar{a}f(b))=0, \]
	where the last equality holds because $ D $ is an intrinsic extension of $ \bar{a}b $. But, the above inequality would contradict the fact that $ \bar{a}E'_0f(E_1) $ is $ n $-weakly closed in $ M_2 $ and hence there would not exist any realization of $ D $ over this embedding.
	
	The discussion about type-2 structures appeared just above shows that in Case B.2, we have to to let $ \Omega $ consist of all type-2 structures $ D'_j $, as defined above, and to let $ n $ be the maximum of all numbers $ \card[D_j\backslash\bar{a}] $ with $ D_j $ ranging over all structures appearing in $\varphi$ (regardless of whether they are type-1 and type-2). Then, an application of the axiom \axsem for $ E'_0, E_1, n $ and $ \Omega $ yields an embedding $ f $ of $ E_1 $ over $ E'_0 $ into $ M_2 $ such that all structures $ D_j $ in $ \varphi(\bar{x}y) $ get omitted over $ f(E_1) $.
	
	To see the completeness of $ \tmu $, note that, by Corollary \ref{crlJEP}, this theory has the joint embedding property, and we just showed that it eliminates quantifiers. For decidability, we can use \Cref{crl-ommitable} to easily verify that the axiom schemes of $T_0,$ \axsem, and \axisem are recursively enumerable. 
	
\end{proof}

\begin{crl}\label{crlModelCompanion}
	$ \tmu $ is a model companion for $ \tplus $. This theory axiomatizes the class of existentially closed models of $ \tplus. $
\end{crl}
\begin{proof}
	By Lemma \ref{lmaExCloedIsUCCGeneric}, any existentially closed model of $ \tplus $ is a model of $ \tmu $. On the other hand, $ \tplus $ is an inductive theory and each of its models can be extended to an existentially closed model. Hence, $ \tmu $ and $ \tplus$ are cotheories. Moreover, $ \tmu $, by Theorem \ref{thmQE}, has quantifier elimination which shows that it is a model companion for $ \tplus $. Since $ \tplus $ is inductive, the latter is equivalent to saying that $ \tmu $ axiomatizes the class of all existentially closed models of $ \tplus $.
\end{proof}

\begin{prop}\label{propModelCompletion}
	$ \tmu $ is the model completion of $ \tplus. $
\end{prop}
\begin{proof}
	By Corollary \ref{crlModelCompanion}, $ \tmu $ is a model companion for $ \tplus $ and we only need to show that $ \tplus $ has the amalgamation property. Let $ M_0, M_1, M_2 $ be models of $ \tplus $ with $ M_0\subseteq_{+} M_1 $ and $ M_0\subseteq_{+} M_2 $. If $ M $ denotes the free join of $ M_1 $ and $ M_2 $ over $ M_0 $, by part (ii) of Lemma \ref{lmaAmalgam}, the structure $ M $ would be a model of $ T_{\nat} $ with $ M_1\subseteq_{+}M $ and $ M_2\subseteq_{+} M $. By Lemma \ref{lmaUCCExtension}, the structure $ \uccext $ is a model of $ \tmuinfinity $ with $ M_1\subseteq_{+} \uccext $ and $ M_2\subseteq_{+}\uccext $. Now, it is easy to expand $ \uccext $ to a model of $ T_{\nat} $, and this completes the proof.
\end{proof}

\vspace*{.3cm}
\section{Properties of $\tmu$}\label{secSOP}


We first prove Lemma \ref{lmaIndisc} which underlines a key property of indiscernible sequences in the models of $ \tmu $.

\begin{assumption}
	We fix a monster model $ \mathcal{\mathcal{M}} $ for the theory $ \tmu $ and will omit all subscripts referring to the ambient model; therefore, instead of notations like $ \chi^{*}_{\mathcal{M}}(B/A) $ or $ \cl_{\mathcal{M}}(A) $ we only write $ \chi^{*}(B/A) $ and $ \cl(A) $. 
\end{assumption}

\begin{lem}\label{lmaIndisc}
	Suppose that $ \{\bar{a}_{i}\}_{i\in\omega} $ is an $ A $-indiscernible sequence. Then,
	\begin{itemize}
		\item[(i)]
		The tuples $ \bar{a}_{i} $ are mutually in free amalgamation over $ A $.
		
		\item[(ii)] 
		The sequence $ \big\{\wcl(\bar{a}_{i})\big\}_{i\in\omega} $ is an indiscernible sequence of finite tuples over $ \wcl(A) $ which are mutually in free amalgamation over $ \wcl(A). $
		
		\item[(iii)]
		If $ \bar{A} $ denotes the weak closure of $ A, $ then for any $ n\in\omega $ the weak closure of the structure $ A\bar{a}_{0}\cdots\bar{a}_{n} $ is equal to 
		\[ \freejoin{\freejoin{\wcl(\bar{a}_{0})}{\bar{A}}{\cdots}}{\bar{A}}{\wcl(\bar{a}_{n})}. \]
	\end{itemize} 
\end{lem}
\begin{proof}
	(i) Without loss of generality, we may assume that the intersection of tuples $ \bar{a}_{i} $ lies inside $ A $ and they are mutually disjoint over $ A. $ 
	
	Let $ r>0 $ denote the number of relations preventing $ \bar{a}_{0} $ from being freely amalgamated with $ \bar{a}_{1} $ over $ A. $ A particular consequence of  indiscernibility is that all tuples $ \bar{a}_{i} $ have an identical predimension. Hence, an easy induction on $ n $ shows that 
	\[ \delta(\bar{a}_{0}\cdots\bar{a}_{n-1}) =n\delta(\bar{a}_{0})-\frac{n(n-1)}{2}r. \]
	Therefore, for every number $ n> \frac{2\delta(\bar{a}_{0})}{r}+2 $ we have that $ \delta(\bar{a}_{0}\cdots\bar{a}_{n-1})<0. $ This would contradict the fact that $ \mathcal{M}\in\kpluszerobar $.
	
	(ii) It is easy to see that $ \big\{\wcl(\bar{a}_{i})\big\}_{i\in\omega} $ is indiscernible over the weak closure of $ A $. The claimed free amalgamation then follows directly from part (i).
	
	(iii) For any $ i\in\omega $ let $ \bar{c}_{i}:= \wcl(\bar{a}_{i}). $ Using parts (i) and (ii), we know that the sequence $ \big\{\bar{c}_{i}\big\}_{i\in\omega} $ is $ \bar{A} $-indiscernible with its elements being in free amalgamation over $ \bar{A} $. Hence, we only need to show that the structure $ \freejoin{\freejoin{\bar{c}_{0}}{\bar{A}}{\cdots}}{\bar{A}}{\bar{c}_{n}} $ is weakly closed. We show this for $ n=2 $ which completely resembles the proof of the general case.
	
	If the structure $ \freejoin{\bar{c}_{0}}{\bar{A}}{\bar{c}_{1}} $ is not weakly closed, then there will exist a finite structure $ B $ and a finite substructure $ A_{0}\subsetfinite\bar{A} $ satisfying
	\[ \delta(B/A_{0}\bar{c}_{0}\bar{c}_{1})<0. \]
	Hence, by indiscernibility, for any $ n\in\omega $ there are $ \frac{n(n-1)}{2} $-many realizations of $ B $ over the structure $ A_{0}\bar{c}_{0}\cdots\bar{c}_{n-1} $. If $ S_{n} $ denotes the union of $ A_{0}\bar{c}_{0}\cdots\bar{c}_{n-1} $ with all these realizations of $ B, $ we will have that
	\begin{align*}
		\delta(S_{n})&= \delta(A_{0}\bar{c}_{0}\cdots\bar{c}_{n-1})+\frac{n(n-1)}{2}\delta(B/A_{0}\bar{c}_{0}\bar{c}_{1})\\
		&\leq n\delta(A_{0}\bar{c}_{0}\bar{c}_{1})+\frac{n(n-1)}{2}\delta(B/A_{0}\bar{c}_{0}\bar{c}_{1}).
	\end{align*} 
	This will lead to a similar contradiction appeared in the proof of part (i) of this lemma.
\end{proof}

\begin{rem}
	
	\Cref{lmaIndisc} seems to 
	be beneficial even in the more general context of an arbitrary amalgamation class built upon a predimension function. In fact, comparing to the subtleties that appear in deploying the dimension theory in such a context, working with this lemma, in the form presented here, seems to be more efficient when one concerns with investigating the existence/non-existence of combinatorial configurations that appear in different model-theoretic properties: Order property, independence property, $\trp, \trp_1, \trp_2, \sop_{n}$, or the strict order property.
	
%
	
\end{rem}

\begin{notation}
	Let $ \tr(x,y)$ be the formula  $\exists z R(x,y,z)$ which asserts that $x$ and $y$ participate in a triangle.
	
\end{notation}

\begin{thm}\label{thm-sop}
	$ \tmu$ has the strict order property.
\end{thm}
\begin{proof}
	
	Let 
	$\varphi(y_1z_1w_1;y_2z_2w_2)$ be the following formula:
	\begin{align}\label{eq-sop}
		\begin{cases}
			y_{1}\neq y_{2}\ \wedge z_1 = z_2 \ \wedge w_1 = w_2\  \wedge \\
			\exists x \big{(}R(x,z_1,w_1)\wedge \tr(x,y_{1}) \wedge\neg\tr(x,y_{2})\big{)} \ \wedge \\
			\forall x \big{(}(R(x,z_1,w_1) \wedge \neg\tr(x,y_{1}))\longrightarrow \neg\tr(x,y_{2}) \big{)}.
		\end{cases}
	\end{align}
	For better readability of the argument, let $\theta(y_1, y_2, z_1,w_1)$ and $\psi(y_1, y_2, z_1,w_1)$ respectively denote the second and the third line of 
	\eqref{eq-sop}. 
	
	It is obvious that $\varphi$ is irreflexive and antisymmetric. To see that $\varphi$ has also the transitive property, suppose that we have $\models \varphi(b_1c_1d_1;b_2c_2d_2)\wedge\varphi(b_2c_2d_2;b_3c_3d_3)$ for some elements of $\mathcal{M}$. First note that the first line of \eqref{eq-sop} implies that $c_1=c_2=c_3$ and $d_1=d_2=d_3$; let $c$ and $d$ denote theses elements respectively. By $\theta(b_1,b_2,c,d)$ and $\theta(b_2,b_3,c,d)$, there are elements $a_1$ and $a_2$ such that $\models R(a_1,c,d)\wedge R(a_2,c,d)$, and $\models \tr(a_1,b_1)\wedge \tr(a_2,b_2)$. It suffices to show that $\models \neg\tr(a_1,b_3)\wedge  \psi(b_1,b_3,c,d)$. For the first formula, note that if $\models \tr(a_1,b_3)$, then, by $\psi(b_2,b_3,c,d)$, we have that $\models \tr(a_1, b_2)$ which contradicts $\psi(b_1,b_2,c,d)$. To see $\models \psi(b_1,b_3,c,d)$, suppose that $\models R(a,c,d)\wedge \neg \tr(a,b_1)$ for some element $a\in\mathcal{M}$. Then, by $\psi(b_1,b_2,c,d)$, we have that $\models \neg \tr(a,b_2)$. This, using $\psi(b_2,b_3,c,d)$, implies that $\models \neg \tr(a,b_3)$. 
	
	Therefore, $\varphi$ defines a partial strict order over $\mathcal{M}$. We only need to show the existence of an infinite sequence 
	$(b_i)_{i<\omega}$ and two elements $c$ and $d$ such that $\models  \varphi(b_icd;b_jcd)$ for all $0\leq i<j<\omega$.
	
	We show that every finite part of this sequence is realizable in $\mathcal{M}$. To do this, for each $n$, we construct a finite structure $A_n\in\kpluszero$, and suitably embed it into $\mathcal{M}$. Let $A_n$ be the structure that consists of the elements 
	\begin{align*}
		&\{b_i: 0\leq i\leq n\}\cup\{c,d\}\cup\\
		&\{a_{i}: 0\leq i< n\}\cup\{d_{i}:0\leq i<n\}\cup\\
		&\{d'_{(i,j)}:0\leq i<j< n\},	
	\end{align*}
	
	in which only the following relations hold:
	\[ \bigwedge_{0\leq i< n}R(c,d,a_{i}) \ \wedge \ \bigwedge_{0\leq i< n} R(b_i,a_{i}, d_{i}) \ \wedge \ \bigwedge_{0\leq i<j< n} R(b_i,a_{j},d'_{(i,j)}). \]
	
	For each $0\leq i< n$, the relations $R(b_i,a_i, d_i)$ ensure the satisfaction of $\theta(b_i,b_j,c,d)$ in $A_n$, as they imply $\tr(b_i,a_i)\wedge \neg\tr(b_j,a_i)$. Moreover, for each $b_i$ and $b_j$ with $0\leq i<j<n$, only the elements $a_j$ are capable of violating $\psi(b_i,b_j,c_1,c_2)$, as there is no other element forming a triangle with $b_j$ in $A_n$. However, the elements $d'_{(i,j)}$ with their corresponding relations ensure that $\tr(b_i,a_j)$ holds as well. This implies that $A_n\models \psi(b_i,b_j,c,d)$.

	Now, for some $i$ consider the pair of structures $(cd,cda_i) $ for which we have $\intrext[cd]{cda_0}$. The number of copies of $a_i$ over $cd$ grows in $A_n$ as we increase $n$. Hence, for all $n>\mu(2,3)$, when we embed $A_n$ in a paracollapsed structure, say $M$, there must exist infinitely many copies of $a_i$ over $cd$ in $M$. But, it is easy to see that these new triangles do not violate the satisfaction of $\varphi(b_icd;b_jcd)$ for all $0\leq i<j<n$ in $M$. Moreover, observe that the pairs $(cd,cda_i)$ alone need such a treatment, as for all other intrinsic extensions available in $A_n$, say  $\intrext[B]{C}$, the number of copies of $C$ over $B$ is equal to 1.
		
	Finally, if we embed $A_n$ in $\mathcal{M}$ sufficiently closed so that no new triangle--other than the ones discussed above--arises over $A_n$ in $\mathcal{M}$ we can conclude that $\mathcal{M}\models \varphi(b_icd;b_jcd)$ for each $0\leq i<j<n$. But, this can be done using the axiom of semigenericity. 

\end{proof}
\begin{dfn}\label{dfn-tp2}
	Let $ T $ be a complete first order theory. A formula $ \varphi(\bar{x};\bar{y}) $ in the language of $T$ has the \textit{tree property of the second kind} ($ \trp_2 $) if there is a model $ M\models T $ and a set of tuples $ \big\{\bar{b}_{i,j}:i,j<\omega\big\}\subseteq M $ such that 
	\begin{itemize}
		\item[(i)] for every $ \sigma\in\omega^\omega $ the set of formulas 
		$ \big\{\varphi(\bar{x};\bar{b}_{n,\sigma(n)}):n\in\omega\big\} $ is consistent, and
		\item[(ii)] for all $n\in\omega $ and $n<i<j<\omega$ the set $ \varphi(\bar{x};\bar{b}_{n, i})\wedge\varphi(\bar{x};\bar{b}_{n, j}) $ is inconsistent.
	\end{itemize}
	$T$ is $\ntp_2$ if no formula can witness $\trp_2$ in models of $T$.
\end{dfn}
More intuitively, for each $n\in\omega $ list all the formulas $\{\varphi(\bar{x};\bar{b}_{n,j}):n<j<\omega\}$ in a row. Then, if we choose exactly one formula out of each row, we gain a consistent set of formulas, while, no two formulas on the same row are consistent. 
\begin{thm}\label{thm-tp2}
	$ \tmu $ has $\trp_2$.
\end{thm}
\begin{proof}
	
	Let $ \varphi(x;y_1y_2y_3) $ be the following formula
	\[ R(x,y_1,y_2)\wedge\tr(x,y_{3}). \]
	
	To show that $ \varphi $ witnesses $ \trp_2 $, we construct a finite structure $ A $  and embed it suitably into $ \mathcal{M} $.

	Fix a tuple $ \bar{b}=b^0b^1 $ and let $ \big{(}\bar{b}_{i,j}\big{)}_{ i,j<\omega } $ be a finite set of tuples, each of size three, with
	\[ \bigcap_{i,j<\omega}\bar{b}_{i,j}=\bar{b}. \]
	For each $ i,j<\omega$, define $ b_{i,j} $ to be the unique element in $ \bar{b}_{ij}\backslash\bar{b} $; that is, $\{b_{i,j}\}$ is the singleton obtained by removing the fixed tuple $\bar{b}$ from $\bar{b}_{i,j}$.
	
	For each $ \sigma\in \omega^{\omega}, $ let $ A_\sigma $ be a structure containing the fixed tuple  $ \bar{b} $ along with a new element  $ a_\sigma$ satisfying $ R(a_\sigma,b^0,b^1) $. Additionally, for each $ n<\omega$, define $ C_{\sigma(n)} $ as the structure extending $ A_\sigma b_{n,\sigma(n)} $ by adding a new element $ c_{\sigma(n)} $ such that $ R(a_\sigma,b_{n,\sigma(n)},c_{\sigma(n)}) $ holds. The element $a_{\sigma}$ provides the branch consistency condition.

	Finally, let $ A $ denote the following structure
	\[ \bigcup_{\sigma\in \omega^\omega}\bigcup_{n<\omega}C_{\sigma(n)}. \]
	
	It is straightforward to verify that $ A $ belongs to $ \kpluszero $.
	
	Now, any structure $ C_{i,j}\in\kpluszero $ that realizes the set of formulas  $\{ \varphi(x;\bar{b}_{n,i}),\varphi(x;\bar{b}_{n,j})\} $ must contain the tuples $ \bar{b}_{n,i}\bar{b}_{n,j} $ along with elements $ a $, $c_1 $ and $ c_2 $ satisfying  
	$ R(a,b^0,b^1)\wedge R(a,b_{n,i},c_1)\wedge R(a,b_{n,j},c_2) $. However, it is clear that such a structure $ C_{i,j} $ is omittable over $ A $, since for any $n\in\omega$ and any $n<i<j<\omega$ the structure $ A $ does not contain even a single copy of $ C_{i,j} $ over $ \bar{b}_{n,i}\bar{b}_{n,j} $. Hence, it suffices to apply \axsem to embed $ A $ into $ \mathcal{M} $ in such a way that it omits all structures $C_{i,j}$. This will guarantee the inconsistency condition required by $\trp_2$.
\end{proof}

\subsection{Formulas Witnessing Other Order Properties}
For each $n\geq 2$, having $\sop_{n}$ for a theory is a consequence of having the strict order property. Although the examples below are not going to prove a new property for $\tmu$, it seems model-theoretically interesting to see explicit examples of formulas that witness different order properties in the context of hypergraphs. 

We first introduce two auxiliary formulas that will be used in Examples \ref{ex-sop3} and \ref{ex-sopn}. Let $\fan(x,y,w_1,w_2,x')$ be the formula $ \tr(x,y)\wedge R(y,w_1,w_2)\wedge \tr(y,x')$ which asserts the existence of three triangles meeting one another at $y$, which form a shape similar to a fan! Also, let $\spath(x,w_1,w_2,x')$ be the formula asserting the existence of two distinct elements $y_1,y_2$ satisfying the following formula:
\begin{align}\label{form-spath}
	\neg\fan(x,y_1,w_1,w_2,x')\wedge \neg\fan(x,y_2,w_1,w_2,x')\ \wedge
	\tr(x,y_1)\wedge \tr(x',y_2).
\end{align}
The latter formula will appear in \Cref{ex-sopn} where we also enforce $R(y_1,w_1,w_2)$ and $R(y_2, w_1,w_2)$ to hold. In that setting, and intuitively speaking, this formula is intended to describe a special kind of path that connects $x$ to $x'$ through $y_1, w_1, w_2 $ and $y_2 $ while preventing the latter elements from  participating in fan-like structures involving $x$ and $x'$.

The first example in due to Alex Kruckman who kindly let us include it in the paper:
\begin{example}\label{ex-sop3}
	A formula witnessing $\sop_3$.
\end{example}
\begin{proof}
	Let $\psi(x,y,w_{1},w_{2})$ be the formula: 
	\begin{align*}
		R(y,w_{1},w_{2})\wedge\tr(x,y).
	\end{align*}
	Let $\varphi(x_{1}, x_{2},w_{1},w_{2})$ be the formula: 
	\begin{align*}
		\exists y \fan( x_{1}, y, w_{1}, w_{2},x_{2})
	\end{align*}
	Writing $\bar{x}$ for the tuple of variables $(x_{1}, x_{2}, y, w_{1}, w_{2})$, let $\theta(\bar{x}, \bar{x}^{'})$ be the formula:
	\begin{align*}
		(w_{1} = w^{'}_{1}) \wedge(w_{2} =w^{'}_{2})\wedge\neg\varphi (x^{'}_{1}, x_{2}, w_{1}, w_{2})\wedge\psi (x_{1}, y^{'}, w_{1}, w_{2}) \wedge\psi(x^{'}_{2}, y, w_{1}, w_{2}) 
	\end{align*}
	We claim that $\theta(\bar{x}, \bar{x}^{'})$ has SOP$_3$. For this, we first need to show that 
	\begin{align*}
		\theta(\bar{x}, \bar{x}^{'})\wedge \theta(\bar{x}^{'}, \bar{x}^{''})\wedge \theta(\bar{x}^{''}, \bar{x})
	\end{align*}
	is inconsistent. Indeed, $\theta(\bar{x}, \bar{x}^{'})$ implies $\neg\varphi (x^{'}_{1}, x_{2}, w_{1}, w_{2})$, but the conjunction above also implies
	\begin{align*}
		\psi(x^{'}_{1}, y^{''}, w_{1}, w_{2}) \wedge \psi(x_{2}, y^{''}, w_{1}, w_{2})
	\end{align*}
	implying $\varphi(x^{'}_{1}, x_{2},w_{1},w_{2})$, which is a contradiction. 
	
	Now it suffices by compactness to find, for each $n\in \omega$ a sequence 
	$(\bar{a}^{i})_{i<n}$ such that $\mathcal{M} \models  \theta(\bar{a}^{i}, \bar{a}^{j})$ for all $0\leq i<j<n$. To do this, we describe a structure $A$ in $\mathcal{K}_{0}$ and use semigenericity to find an appropriate embedding of $A$ into $\mathcal{M}$.
	
	The structure $A$ has distinct elements $(a^{i}_{1})_{i<n}, (a^{i}_{2})_{i<n}, (b^{i})_{i<n}, (c^{i,j}_{1})_{i<j<n}, (c^{i,j}_{2})_{i<j<n}, d_{1}$ and $d_{2}$ for a total of $3n + n(n 1)+2$ elements. The relation $R^{A}$ consists of the following triples:
	\begin{align*}
		&\{b^{i},d_{1},d_{2}\} \hspace*{1mm}\text{for all} \hspace*{1mm} i<n\\
		&\{a^{i}_{1},b^{j},c^{i,j}_{1}\} \hspace*{1mm}\text{for all} \hspace*{1mm} i<j<n\\
		&\{a^{j}_{2},b^{i},c^{i,j}_{2}\} \hspace*{1mm}\text{for all} \hspace*{1mm} i<j<n
	\end{align*}
	To check that $A \in \mathcal{K}_{0}$, note that each triple of the form $\{a^{i}_{1}, b^{j}, c^{i,j}_{1}\}$ or $ \{a^{j}_{2}, b^{i},c^{i,j}_{2}\}$ contains an element which is not contained in any other triple, so triples of this form cannot contribute to a non-positive predimension. Similarly, each triple of the form $\{b^{i}, d_{1}, d_{2}\}$ contains an element which is not contained in any other triple of this form.
	
	We define $d^{i}_{1} = d_{1}$ and $d^{i}_{2} = d_{2}$, and write $\bar{a}^{i} = (a^{i}_{1}, a^{i}_{2},  b_{i}, d^{i}_{1}, d^{i}_{2})$, for all $i<n$. Then for all $i <j<n$, we have $A\models \theta(\bar{a}^{i},\bar{a}^{j})$.  With the exception of $\neg\varphi(x^{'}_{1},x_{2},w_{1},w_{2})$  each of the other four conjuncts of $\theta $ are existential, so their satisfaction is preserved by any embedding of $A$ into $\mathcal{M}$. It remains to show that we can find some embedding such that $\mathcal{M} \models \neg\varphi (a^{j}_{1}, a^{i}_{2}, d_{1}, d_{2})$ for all $i<j<n$.
	
	Fix $i<j<n$ and let $C$ be any extension of $A$ by elements called $b^{\ast}$, $c^{\ast}_{1}$, and   $c^{\ast}_{2}$ (we allow $c^{\ast}_{1}=c^{\ast}_{2}$, or for some of these elements to be equal to elements already in $A$), and any number of triples involving the new elements, at least including $\{b^{\ast}, d_{1}, d_{2}\}$, $\{a^{j}_{1}, b^{\ast}, c^{\ast}_{1}\}$, and $\{a^{i}_{2}, b^{\ast}, c^{\ast}_{2}\}$ . Note that if $\mathcal{M} \models \varphi(a^{j}_{1}, a^{i}_{2}, d_{1}, d_{2})$, this is witnessed by some structure $C$ of this form over $A$. Let $C^{'} = cl_{C}(A)$. Since $\delta(C) \leq \delta(A)$, $C^{'}$ is a proper intrinsic extension of $A$ contained in $C$. 
	
	Let $\Omega$ be the set of all such structures $C^{'}$, for all the finitely many choices of $i <j<n$ and structure $C$. By semigenericity, we can embed $A$ into $\mathcal{M}$ in such a way that $\mathcal{M}$ omits all the structures in $\Omega$ over $A$. In particular, $\mathcal{M}$ omits all the structures $C$ enumerated above over $A$. And this means exactly that $\mathcal{M} \models \neg \varphi(a^{i}_{1}, a^{i}_{2}, d_{1},d_{2})$ for all $i<j<n$ as was to be shown.

\end{proof}

The following example is interesting as it provides a formula that witnesses all strong order properties (for $n\geq 3$) at once:
\begin{example} \label{ex-sopn}
	A formula witnessing $\sop_{n}$ for all $n\geq 3$.
\end{example}
\begin{proof}
	Let $\bar{x}$ abbreviate the tuple  $xyw_1w_2$, and $\varphi(\bar{x},\bar{x}')$ denote the following formula:
	\begin{align}\label{form-phi}
		\begin{cases}
			w_1=w'_1 \wedge w_2=w'_2 \ \wedge\\
			R(y,w_1,w_2)\wedge R(y',w_1,w_2) \wedge\tr(x,y')\ \wedge \\
			\text{ ``All other variables are distinct''}\ \wedge \\
			\text{ ``No other relations exist among the elements appearing in } \bar{x}\bar{x}'\text{''} \ \wedge\\
			\neg\tr(x,y)\wedge \neg\tr(x',y') \ \wedge\\
			\neg\spath(x,w_1,w_2,x').
		\end{cases}
	\end{align}
	
	Note that $\neg\tr(x',y)$ is a logical consequence of $\varphi(\bar{x},\bar{x}')$.
	
	We claim that $\varphi(\bar{x},\bar{x}')$ has $\sop_n$ for every $n\geq 3$. Fix $n\in\omega$ and let $d_1, d_2$ denote two elements. For each  $0<i\leq n$ let $\bar{a}_i$ denote the tuple consisting of the elements $a_ib_id_1d_2$. For each $0<i<j\leq n$ add the fresh elements $c_{ij}$. Let $A_n$ be the structure consisting of the elements $\{\bar{a}_i\}_{1\leq i\leq n}\cup\{c_{ij}\}_{1\leq i<j\leq n}$ satisfying:
	\[ \bigwedge_{0< i\leq n} R(b_i,d_1,d_2)\wedge \bigwedge_{0< i<j\leq n} R(a_i,c_{ij},b_j). \]
	It is easy to verify that $A_n$ belongs to $\kpluszero$ as a finite structure. It is also obvious that for each $0<i<j\leq n$, we have that 
	\[ A_n\models \neg\tr(a_i,b_i)\wedge \neg\tr(a_j,b_i) . \]
	Therefore, provided that there are no elements within $A_n$ witnessing $\spath(a_i,d_1,d_2,a_j)$ for some  $0<i<j\leq n$, we can use \axsem to embed $A_n$ into $\mathcal{M}$ sufficiently closed so as to ensure the satisfaction of $\neg\tr(a_i,b_i)\wedge \neg\spath(a_i,d_1,d_2,a_j)$ in $\mathcal{M}$ for each $0<i<j\leq n$.
	
	Hence, we only need to show that $A_n$ satisfies $\neg\spath(a_i,d_1,d_2,a_j)$ for all positive $i<j\leq n$. To see this for such $i<j$, note that $A_n\models\tr(a_i,b_k)\wedge\tr(a_j,b_k)$ whenever $k>j$, and $A_n\models \neg\tr(a_j,b_k)$ whenever $k\leq j$. That is, there is no way in $A_n$ to construct the special path, asserted by $\spath(a_i,d_1,d_2,a_j)$, between $a_i$ and $a_j$.
	
	By compactness, this proves that there is an infinite chain $\{\bar{a}_i\}_{i\in\omega}$ in $\mathcal{M}$ satisfying $\varphi(\bar{a}_i,\bar{a}_j)$ for each $ i<j\in\omega$.
	
	Now, suppose that there is a finite sequence $\{\bar{a}_i\}_{0< i\leq n}$ in $\mathcal{M}$ witnessing the following cycle of length $n$:
	\[ \mathcal{M}\models \varphi(\bar{a}_1,\bar{a}_2) \wedge \cdots\wedge \varphi(\bar{a}_{n-1},\bar{a}_n) \wedge \varphi(\bar{a}_n,\bar{a}_1).\]
	
	This would imply that $\tr(a_{n-1},b_n)\wedge \tr(a_n,b_1)$ which witnesses $\spath(a_{n-1},d_1,d_2,a_n)$ by letting $y_1=b_n$ and $y_2=b_1$ in Formula \eqref{form-spath}. This would contradict $\varphi(\bar{a}_{n-1},\bar{a}_n)$.

	
\end{proof}

\begin{rem}\label{rem-sopn-nsop}
	Note in particular that the formula $\varphi(\bar{x},\bar{x}')$ introduced in \Cref{ex-sopn} does not globally define a partial order. In fact, $
	\varphi(\bar{x}_1,\bar{x}_2) \wedge \varphi(\bar{x}_2,\bar{x}_3) $ implies the positive part of $\varphi(\bar{x}_1,\bar{x}_3)$, namely the first two lines of \eqref{form-phi}, but not necessarily the remaining (negative) parts of $\varphi(\bar{x}_1,\bar{x}_3)$. Hence, this formula does not readily witness the strict order property. 
\end{rem}

\begin{rem}\hfill
	
	\begin{itemize}
		\item[(1)] 
		It is shown in the proof of Theorem 5.18 in \cite{Pourmahd-SimpleGen} that the formula $\exists z R(x,y,z)$ has the independence property. The same proof works in the setting of $\tmu$.
		
		\item[(2)] Using a similar argument and by making the required adjustments, one can show that the formula introduced in \Cref{thm-tp2} witnesses $\sop_2$ as well. 
	\end{itemize}
	
\end{rem}

\vspace*{.3cm}
\textbf{Acknowledgment.} We would like to thank Alex Kruckman, Nick Ramsey, Gabriel Conant, and Omer Mermelstein for their careful treatment of an earlier version of this paper. In particular, we would like to thank Alex Kruckman who noticed us of a gap in our proofs, showing that $ \tmu $ is not  $ \nsop_3 $. His counter example (\Cref{ex-sop3}) was important in finding the formula introduced in \Cref{ex-sopn}. Part of this work was developed during the time the second author was a visiting research fellow at the School of Mathematics in IPM.

\end{document}